\documentclass[10pt,a4paper,reqno]{amsart}

\usepackage{todonotes}
\usepackage{hyperref}
\hypersetup{
	linktocpage=false,
	pdfauthor={Steven Charlton},
	pdftitle={The alternating block decomposition of iterated integrals, and cyclic insertion on multiple zeta values}		
}
\ifpdf\else
\usepackage[hyphenbreaks]{breakurl}
\fi

\newif\ifprintextra
\printextrafalse

\iftrue
\makeatletter
\def\@settitle{%
  \vspace*{-2em}
  \begin{flushleft}%
    \LARGE\bfseries
    \strut\@title\strut
  \end{flushleft}%
}
\def\@setauthors{%
  \begingroup
  \def\thanks{\protect\thanks@warning}%
  \trivlist
  \raggedright
  \large \@topsep27\p@\relax
  \advance\@topsep by -\baselineskip
  \item\relax
  \author@andify\authors
  \def\\{\protect\linebreak}%
  \authors
  \ifx\@empty\contribs
  \else
    ,\penalty-3 \space \@setcontribs
    \@closetoccontribs
  \fi
  \normalfont
  \endtrivlist
  \endgroup
}
\def\@setaddresses{\par
  \nobreak \begingroup
  \small\raggedright
  \def\author##1{\nobreak\addvspace\smallskipamount}%
  \def\\{\unskip, \ignorespaces}%
  \interlinepenalty\@M
  \def\address##1##2{\begingroup
    \par\addvspace\bigskipamount\noindent
    \@ifnotempty{##1}{(\ignorespaces##1\unskip) }%
    {\ignorespaces##2}\par\endgroup}%
  \def\curraddr##1##2{\begingroup
    \@ifnotempty{##2}{\nobreak\noindent\curraddrname
      \@ifnotempty{##1}{, \ignorespaces##1\unskip}\/:\space
      ##2\par}\endgroup}%
  \def\email##1##2{\begingroup
    \@ifnotempty{##2}{\nobreak\noindent E-mail address%
      \@ifnotempty{##1}{, \ignorespaces##1\unskip}\/:\space
      \ttfamily##2\par}\endgroup}%
  \def\urladdr##1##2{\begingroup
    \def~{\char`\~}%
    \@ifnotempty{##2}{\nobreak\noindent\urladdrname
      \@ifnotempty{##1}{, \ignorespaces##1\unskip}\/:\space
      \ttfamily##2\par}\endgroup}%
  \addresses
  \endgroup
  \global\let\addresses=\@empty
}
\def\@setabstracta{%
    \ifvoid\abstractbox
  \else
    \skip@17pt \advance\skip@-\lastskip
    \advance\skip@-\baselineskip \vskip\skip@
    \box\abstractbox
    \prevdepth\z@ 
    \vskip-28pt
  \fi
}
\renewenvironment{abstract}{%
  \ifx\maketitle\relax
    \ClassWarning{\@classname}{Abstract should precede
      \protect\maketitle\space in AMS document classes; reported}%
  \fi
  \global\setbox\abstractbox=\vtop \bgroup
    \normalfont\small
    \list{}{\labelwidth\z@
      \leftmargin0pc \rightmargin\leftmargin
      \listparindent\normalparindent \itemindent\z@
      \parsep\z@ \@plus\p@
      
    }%
    \item[\hskip\labelsep\bfseries\abstractname.]%
}{%
  \endlist\egroup
  \ifx\@setabstract\relax \@setabstracta \fi
}

\def\ps@headings{\ps@empty
  \def\@evenhead{%
    \setTrue{runhead}%
    \normalfont\scriptsize
    \rlap{\thepage}\hfill
    \def\thanks{\protect\thanks@warning}%
    \leftmark{}{}}%
  \def\@oddhead{%
    \setTrue{runhead}%
    \normalfont\scriptsize
    \def\thanks{\protect\thanks@warning}%
    \rightmark{}{}\hfill \llap{\thepage}}%
  \let\@mkboth\markboth
}\ps@headings

\def\section{\@startsection{section}{1}%
  \z@{-1.4\linespacing\@plus-.5\linespacing}{.8\linespacing}%
  {\normalfont\bfseries\Large}}
\def\subsection{\@startsection{subsection}{2}%
  \z@{-.8\linespacing\@plus-.3\linespacing}{.5\linespacing\@plus.2\linespacing}%
  {\normalfont\bfseries\large}}
\def\subsubsection{\@startsection{subsubsection}{3}%
  \z@{.7\linespacing\@plus.2\linespacing}{-1.5ex}%
  {\normalfont\bfseries}}
\def\paragraph{\@startsection{paragraph}{4}%
  \z@{.7\linespacing\@plus.2\linespacing}{-1.5ex}%
  {\normalfont\itshape}}
\def\@secnumfont{\bfseries}

\renewcommand\contentsnamefont{\bfseries}
\def\@starttoc#1#2{\begingroup
  \setTrue{#1}%
  \par\removelastskip\vskip\z@skip
  \@startsection{}\@M\z@{\linespacing\@plus\linespacing}%
    {.5\linespacing}{
      \contentsnamefont}{#2}%
  \ifx\contentsname#2%
  \else \addcontentsline{toc}{section}{#2}\fi
  \makeatletter
  \@input{\jobname.#1}%
  \if@filesw
    \@xp\newwrite\csname tf@#1\endcsname
    \immediate\@xp\openout\csname tf@#1\endcsname \jobname.#1\relax
  \fi
  \global\@nobreakfalse \endgroup
  \addvspace{32\p@\@plus14\p@}%
  \let\tableofcontents\relax
}
\def\contentsname{Contents}
\def\l@section{\@tocline{2}{.5ex}{0mm}{5pc}{}}
\def\l@subsection{\@tocline{2}{0pt}{2em}{5pc}{}}
\makeatother
\fi 
\usepackage{etoolbox}
\makeatletter
\expandafter\patchcmd\csname\string\proof\endcsname{\itshape}{\bfseries}{}{}

\usepackage[marginratio=1:1]{geometry}

\usepackage[backend=bibtex,maxbibnames=10,maxcitenames=5,bibstyle=alphabetic,citestyle=alphabetic]{biblatex}
\defbibheading{bibliography}[\refname]{%
	\section*{#1}%
}

\usepackage{lmodern}
\usepackage[T1]{fontenc}

\usepackage{multirow}
\usepackage{graphicx}
\usepackage{xcolor}

\newcommand{\mailto}[1]{\href{mailto:#1}{#1}}
\newcommand\numberthis{\addtocounter{equation}{1}\tag{\theequation}}

\usepackage{enumitem}
\setlist[itemize]{itemsep=0.5em}
\setlist[enumerate]{itemsep=0.5em}
\setlist[enumerate,1]{label=\roman*)}

\usepackage{mathtools}
\usepackage{thmtools}

\declaretheorem[
	style=plain,
	name=Theorem,
	numberwithin=section,
	refname={Theorem,Theorems},
	Refname={Theorem,Theorems}
]{Thm}
\declaretheorem[
	style=plain,
	name=Proposition,
	numberlike=Thm,
	refname={Proposition,Propositions},
	Refname={Proposition,Propositions}
]{Prop}
\declaretheorem[
	style=definition,
	name=Property,
	numberlike=Thm,
	refname={Property,Properties},
	Refname={Property,Properties}
]{Ppty}
\declaretheorem[
	style=plain,
	name=Corollary,
	numberlike=Thm,
	refname={Corollary,Corollaries},
	Refname={Corollary,Corollaries}
]{Cor}
\declaretheorem[
	style=plain,
	name=Conjecture,
	numberlike=Thm,
	refname={Conjecture,Conjectures},
	Refname={Conjecture,Conjectures}
]{Conj}
\declaretheorem[
	style=plain,
	name=Lemma,
	numberlike=Thm,
	refname={Lemma,Lemmas},
	Refname={Lemma,Lemmas}
]{Lem}
\declaretheorem[
	style=definition,
	name=Definition,
	numberlike=Thm,
	refname={Definition,Definitions},
	Refname={Definition,Definitions}
]{Def}
\declaretheorem[
	style=definition,
	name=Notation,
	numberlike=Thm,
	refname={Notation,Notations},
	Refname={Notation,Notations}
]{Not}
\declaretheorem[
	style=definition,
	name=Example,
	numberlike=Thm,
	refname={Example,Examples},
	Refname={Example,Examples}
]{Eg}
\declaretheorem[
	style=definition,
	name=Remark,
	numberlike=Thm,
	refname={Remark,Remarks},
	Refname={Remark,Remarks}
]{Rem}
\declaretheorem[
	style=definition,
	name=Observation,
	numberlike=Thm,
	refname={Observation,Observations},
	Refname={Observation,Observations}
]{Obs}

\allowdisplaybreaks

\usepackage{braket}
\usepackage{shuffle}

\renewcommand{\vec}[1]{\underline{\mathbf{#1}}}
\newcommand{\Q}{\mathbb{Q}}

\newcommand{\R}{\mathbb{R}}
\newcommand{\wt}{\mathrm{wt}}

\newcommand{\eqN}{\overset{?}{=}}
\newcommand{\eqQ}{\overset{\Q}{=}}
\newcommand{\eqQone}{\overset{1}{=}}

\DeclareMathOperator{\block}{block}
\DeclareMathOperator{\word}{word}
\newcommand{\bl}{\mathrm{bl}}
\newcommand{\st}{\mathrm{st}}
\newcommand{\en}{\mathrm{en}}
\newcommand{\len}{\mathrm{L}}

\newcommand{\gmot}{\mathfrak{a}}
\newcommand{\mot}{\mathfrak{m}}
\newcommand{\lmot}{\mathfrak{L}}
\newcommand{\dec}{\mathrm{dec}}

\newcommand{\refl}{\operatorname{\mathcal{R}}}

\newcommand{\size}[1]{\#{#1}}

\newcommand{\argdot}{{{}\mathop{\cdot}{}}}

\DeclareMathOperator{\per}{per}
\DeclareMathOperator{\id}{id}

\DeclareMathOperator{\Sym}{Sym}
\DeclareMathOperator{\Alt}{Alt}

\newcommand{\cyc}{\operatorname{\mathcal{C}}}

\renewcommand{\epsilon}{\varepsilon}

\makeatletter
\let\@@pmod\pmod
\DeclareRobustCommand{\pmod}{\@ifstar\@pmods\@@pmod}
\def\@pmods#1{\mkern4mu({\operator@font mod}\mkern 6mu#1)}
\makeatother

\makeatletter
\def\proofstep{\@startsection{paragraph}{4}%
	\z@{1ex}{-\fontdimen2\font}%
	{\normalfont\itshape}}
\makeatother

\begin{document}
	
	\title[Alternating block decomposition, and cyclic insertion on MZV's]{The alternating block decomposition of iterated integrals, \\ and cyclic insertion on multiple zeta values}
	\author{Steven Charlton}
	
	\address{Fachbereich Mathematik \\	
		Auf der Morgenstelle 10 (Geb\"aude C) \\
		72076 \\
		T\"ubingen \\
		Germany}
	
	\email{\mailto{charlton@math.uni-tuebingen.de}}
	
	
	\date{27 April 2017}
	
	\begin{abstract}
		The cyclic insertion conjecture of Borwein, Bradley, Broadhurst and Lison\v{e}k states that by inserting all cyclic permutations of some initial blocks of 2's into the multiple zeta value \( \zeta(1,3,\ldots,1,3) \) and summing, one obtains an explicit rational multiple of a power of \( \pi \).  Hoffman gives a conjectural identity of a similar flavour concerning \( 2 \zeta(3,3,\{2\}^m) - \zeta(3,\{2\}^m,(1,2)) \).
				
		In this paper we introduce the `\emph{generalised} cyclic insertion conjecture', which we describe using a new combinatorial structure on iterated integrals -- the so-called \emph{alternating block} decomposition.  We see that both the original BBBL cyclic insertion conjecture, and Hoffman's conjectural identity, are special cases of this \emph{generalised} cyclic insertion conjecture.  By using Brown's motivic MZV framework, we establish that some symmetrised version of the generalised cyclic insertion conjecture always holds, up to a rational; this provides some evidence for the generalised conjecture.
	\end{abstract}
	
	\maketitle

	\tableofcontents

	\section{Introduction}
	
	\subsection{Previous results on cyclic insertion}
	
	Throughout \cite{borwein1998combinatorial}, \citeauthor{borwein1998combinatorial} present numerical evidence for a `cyclic insertion conjecture' on multiple zeta values, Conjecture 1 in \cite[][p. 9]{borwein1998combinatorial}.  Their conjecture can be given as follows
	
	\begin{Conj}[BBBL cyclic insertion, Conjecture 1 in {\cite{borwein1998combinatorial}}]\label{conj:bbblcyclic}
		For given non-negative integers \( b_0, b_1, \ldots, b_{2n} \geq 0 \), we have
		\[
			\sum_{\sigma \in C_{2n+1}} \zeta(\{2\}^{b_{\sigma(0)}}, 1, \{2\}^{b_{\sigma(1)}}, 3, \ldots, 1, \{2\}^{b_{\sigma(2n-1)}}, 3, \{2\}^{b_{\sigma(2n)}}) \eqN \frac{\pi^\wt}{(\wt+1)!} \, ,
		\]
		where \( C_{2n+1} = \langle (0 \, 1 \, \cdots \, 2n) \rangle \) is the cyclic group of order \( 2n+1 \) viewed as a subgroup of the symmetric group on \( \Set{0,1,\ldots,2n} \).  Here \( \wt \) denotes the weight of the multiple zeta values on the left hand side, which in this case is \( 2 \sum_{i=0}^{2n} b_i + 4 n \).  We also use the shorthand notation \( \{s\}^{\ell} \coloneqq \underbrace{s, s, \ldots, s}_{\text{\( \ell \) times}} \). \medskip
		
		That is, the fixed blocks \( \{2\}^{b_0},\ldots,\{2\}^{b_{2n}} \) are inserted between the arguments of \( \zeta(\{1,3\}^n) \), and when all cyclic shifts of these blocks are summed, the result is \( \frac{\pi^\wt}{(\wt+1)!} \).
	\end{Conj}
	
	Throughout this paper we will abbreviate multiple zeta value to MZV.  We will keep with the convention that
	\[
		\zeta(s_1,\ldots,s_k) \coloneqq \sum_{n_1 < n_2 < \cdots < n_k} \frac{1}{n_1^{s_1} \cdots n_k^{s_k}} \, ,
	\]
	so that \( \zeta(1,2) \) is a \emph{convergent} MZV.
	
	The \emph{weight} of \( \zeta(s_1,\ldots,s_k) \) is the sum \( \sum_i s_i \) of the arguments, and the \emph{depth} of \( \zeta(s_1, \ldots, s_k) \) is the number \( k \) of the arguments.  We will write `\( \wt \)' on the right hand side of identities as shorthand for the weight of the (homogeneous) combination of MZV's on the left hand side.  We also employ the shorthand notation
	\[
		\{ x_1, x_2, \ldots, x_n \}^\ell \coloneqq \underbrace{x_1, x_2, \ldots, x_n, \ldots, x_1, x_2, \ldots, x_n}_{\text{\( \ell \) copies of \( x_1,x_2,\ldots, x_n \)}} \, ,
	\] to write repeated arguments.  \medskip

	In the special case where all \( b_0 = \cdots = b_{2n} = 0 \), the BBBL cyclic insertion conjecture reduces to the Zagier-Broadhurst identity
	\[
		\zeta(\{1,3\}^n) = \frac{1}{2n+1} \frac{\pi^\wt}{(\wt+1)!} \, ,
	\]
	conjectured by Zagier, and first proven by Broadhurst in Section 11 of \cite{borwein2001special} using hypergeometric functions.  In the special case where \( n = 1 \), and \( b_1 = m \), we obtain the known evaluation
	\begin{equation}\label{eq:zeta222}
		\zeta(\{2\}^m) = \frac{\pi^\wt}{(\wt+1)} \, ,
	\end{equation}
	as given in Equation 36 of \cite[][p. 7]{borwein1997evaluations}.  This result is proven by generating function methods, by generalising Euler's evaluation of \( \zeta(2) \).  Nevertheless, the slightly more general case \( b_0 = \cdots = b_{2n} = m \) leads to following conjectural evaluation
	\begin{equation}\label{eqn:zeta21232:eval}
		\zeta(\{\{2\}^m, 1, \{2\}^m, 3\}^n, \{2\}^m) \eqN \frac{1}{2n+1} \frac{\pi^\wt}{(\wt+1)!} \, ,
	\end{equation}
	proposed earlier in Equation 18 of \cite[][p. 4]{borwein1997evaluations}.
	
	The special case \( b_0 = 1 \), and \( b_1 = \cdots = b_{2n} = 0 \) of the BBBL cyclic insertion conjecture is already proven by \citeauthor{borwein1998combinatorial} in Theorem 2 of the original paper \cite[][p. 8]{borwein1998combinatorial}.  They say this gives a version of Zagier's identity ``dressed with 2's''. \medskip
	
	In \cite{bowman2002algebra}, \citeauthor{bowman2002algebra} generalise the result for \( b_0 = 1 \), \( b_1 = \cdots = b_{2n} = 0 \) to a sum over all compositions \( \sum_{i=0}^{2n} b_i = m \), giving the Bowman-Bradley Theorem.
	
	\begin{Thm}[Bowman-Bradley, Corollary 5.1 in \cite{bowman2002algebra}]\label{thm:bowmanbradley}
		Let \( n, m \) be non-negative integers.  Then
		\begin{align*}
			\sum_{\substack{b_0 + \cdots + b_{2n} = m \\ a_i \geq 0 }} \zeta(\{2\}^{b_0}, 1, \{2\}^{b_{1}}, 3, \ldots, 1, \{2\}^{b_{2n-1}}, 3, \{2\}^{b_{2n}}) \hspace{2em} \\ 
			{} = \frac{1}{2n+1} \binom{m + 2n}{m} \frac{\pi^\wt}{(\wt+1)!} \, .
		\end{align*}
		
		That is, blocks \( \{2\}^{b_0}, \ldots, \{2\}^{b_{2n}} \) corresponding to all possible compositions (including those where some \( b_i = 0 \) are permitted) of \( m \) into \( 2n+1 \) parts are inserted into \( \zeta(\{1,3\}^n) \), and the results are summed.
	\end{Thm}
	
	Simpler and more refined proofs of the Bowman-Bradley theorem have since been given by Zhao \cite{zhao2008exotic}, and Muneta \cite{muneta2009evaluations}.

	The Bowman-Bradley is compatible with the cyclic insertion conjecture because any composition \( m = \sum_{i=0}^{2n} b_i \) of \( m \) into \( 2n+1 \) parts remains a composition of \( m \) into \( 2n+1 \) parts when the parts are cyclically shifted.  One also computes that there are \( \binom{m + 2n}{m} \) such compositions.  According to the BBBL cyclic insertion conjecture, each term in the above sum should contribute on average \( \frac{1}{2n+1} \frac{\pi^\wt}{(\wt+1)} \); we therefore obtain the expected coefficient.  \medskip
	
	In \cite{charlton2015rational}, we proved using Brown's motivic MZV framework \cite{brown2012mixed,brown2012decomposition}, that a sum over all permutations of given blocks \( \{2\}^{b_0}, \ldots, \{2\}^{b_{2n+1}} \) is already sufficient to obtain a rational multiple of \( \pi^\wt \).
	
	\begin{Thm}[Symmetric insertion, Theorem 3.1 in \cite{charlton2015rational}]
		For given non-negative integers \( b_0, b_1, \ldots, b_{2n} \geq 0 \), we have
		\[
		\sum_{\sigma \in S_{2n+1}} \zeta(\{2\}^{b_{\sigma(0)}}, 1, \{2\}^{b_{\sigma(1)}}, 3, \ldots, 1, \{2\}^{b_{\sigma(2n-1)}}, 3, \{2\}^{b_{\sigma(2n)}}) \in \pi^\wt \Q \, ,
		\]
		where \( S_{2n+1} \) is the symmetric group on \( \Set{ 0, 1, \ldots, 2n } \).  \medskip
		
		That is, the fixed blocks \( \{2\}^{b_0},\ldots,\{2\}^{b_{2n}} \) are inserted between the arguments of \( \zeta(\{1,3\}^n) \).  When all permutations of these blocks are summed, the result \emph{is} a rational multiple of \( \pi^\wt \).
	\end{Thm}
	
	As a consequence, in Corollary 3.19 of \cite{charlton2015rational}, we obtained that the evaluation in \autoref{eqn:zeta21232:eval} does indeed hold up to a rational.  Namely
	\[
		\zeta(\{\{2\}^m,1,\{2\}^m,3\}^n,\{2\}^n) \in \pi^\wt \Q \, .
	\]\medskip

	In an apparently unrelated direction, we now recall a different conjectural identity.  Equation 5.6 of \cite[][p. 18]{borweinLectures} attributes the following conjectural identity to Hoffman.
	
	\begin{Conj}[Hoffman]\label{conj:hoffman}
		Let \( m \) be a non-negative integer, then we have
		\[
			2 \zeta(3,3,\{2\}^m) - \zeta(3,\{2\}^m,1,2) \eqN -\frac{\pi^\wt}{(\wt+1)!} \, .
		\]
	\end{Conj}
	
	After a little further scrutiny, one begins to see that this identity is somehow similar to, or at least is in much the same spirit as, the BBBL cyclic insertion conjecture.  In this identity, a fixed block \( \{2\}^m \) is inserted into some MZV's, and the resulting sum is (up to sign), \( \frac{\pi^\wt}{(\wt+1)!} \).
	
	According to \cite{borweinLectures}, this result has been checked up to weight 22, where \( m = 8 \), by Vermaseren, using tables of known MZV relations from the MZV datamine \cite{blumlein2010multiple}.   Moreover following Conjecture 5.10.1 in \cite[][p. 188]{zhao2016multiple}, Zhao notes that by a ``tedious computation'' using the idea of \cite{charlton2015rational}, Hoffman's conjecture can be proven up to \( \Q \).  
	
	\begin{Rem}
		Using the properties of iterated integrals on \( \mathbb{P}^1 \setminus \Set{0, 1, z, \infty} \) from \cite{hirose2017dualitysum}, \citeauthor{hirose2017hoffman} \cite{hirose2017hoffman} have since proven \autoref{conj:hoffman}, and indeed the generalisation from \autoref{eg:hoffman33} to 
		\begin{align*}
		\zeta(\{2\}^a, 3, \{2\}^b, 3, \{2\}^c) 
		-\zeta(\{2\}^b, 3, \{2\}^c, 1, 2, \{2\}^a) + {} \\
		{} + \zeta(\{2\}^c, 1, 2, \{2\}^a, 1, 2, \{2\}^b)
		\eqN -\zeta(\{2\}^{3 + a + b + c}) \, ,
		\end{align*}
		exactly.  This now fixes the rational in \autoref{thm:hoffman33} to \( 1 \).
	\end{Rem}
	
	It turns out that the resemblance between the BBBL cyclic insertion conjecture, and Hoffman's conjecture is much more than merely superficial.  Both of these conjectures arise as special cases of the \emph{generalised cyclic insertion conjecture}, \autoref{conj:generalcyclicinsertion} below.  Moreover, as we see already in \autoref{cor:symins}, a symmetrised version of these generalised cyclic insertion identities can always be proven, up to \( \Q \), using Brown's motivic MZV framework; the result in \cite{charlton2015rational} is a special case of this.  This will give us a framework from which to obtain a new proof of Hoffman's identity, up to \( \Q \), as well as proofs of many other such new identities.
	
	\subsection{Main results and structure of the paper}
	
	This paper has two main goals.  As the first goal, we want to elucidate the \emph{alternating block decomposition} framework which is used to unite the BBBL cyclic insertion conjecture and Hoffman's identity into a single \emph{generalised cyclic insertion conjecture}.  For the second goal, we apply Brown's motivic MZV framework to see how identities may be \emph{generated} from the alternating block decomposition; we want to use this to establish that a symmetrised version of the generalised cyclic insertion conjecture always holds, up to \( \Q \). \medskip
	
	In \autoref{sec:motivicmzvbackground} we briefly review the main properties of Goncharov's motivic iterated integrals, and Brown's motivic MZV framework.  In particular, we recall the derivation operators \( D_r \), and Brown's theorem on \( \ker D_{<N} \), which serve as the main tools for our results.  In \autoref{sec:alternatingblock} we then define the alternating block decomposition of an iterated integral.
	
	\subsubsection*{Alternating block decomposition} The alternating block decomposition takes an iterated integral \( I(w) \), where \( w \) is a word over the alphabet \( \Set{0,1} \), and decomposes \( w \) as a concatenation of alternating blocks \( 0101\cdots \) or \( 1010 \cdots \), by cutting at a repeated letter.  By recording the lengths \( \ell_1, \ldots, \ell_n \) of these blocks (and, if necessary, the starting letter \( \epsilon_1 \) that corresponds to the lower bound of integration), we completely capture the information in \( I(w) \).  We write this information as \( I_\bl(\epsilon_1; \ell_1, \ldots, \ell_n)\).  \medskip
	
	In \autoref{sec:reflectionoperators} we introduce the \emph{reflection operators} \( \refl_{j,k} \), the \emph{subsequence} reflection operator \( \refl \) and the notion of a reflectively closed set.  These are employed in \autoref{sec:identitiesfromreflection} as tools to enable the rigorous computation of Brown's derivations \( D_r \) on the sum of iterated integrals of a reflectively closed set.  We prove \( D_r \) always vanishes on such a sum, and so we obtain the main theorem of this paper and the following useful corollary
	
	\begin{restatable*}{Thm}{thmreflectivelyclosed}
		\label{thm:reflclosedid}
		Let \( S \) be a reflectively closed set of block decompositions, with some fixed weight \( N \), and fixed number of blocks \( n \).  Then the sum of the corresponding block integrals satisfies the following
		\[
		\sum_{s \in S} I_\bl^\mot(s) \in \zeta^\mot(N) \Q
		\]
	\end{restatable*}
	
	\begin{restatable*}[Symmetric insertion]{Cor}{corsymmetricinsertion}
		\label{cor:symins}
		Let \( B = (0; \ell_1, \ldots, \ell_n) \) be a block decomposition of (even) weight \( N \geq 2 \).  Then
		\[
		\sum_{\sigma \in S_n} I_\bl^\mot(\ell_{\sigma(1)}, \ldots, \ell_{\sigma(n)}) \in \zeta^\mot(N) \Q
		\]
	\end{restatable*}
	
	In \autoref{sec:cyclicinsertion:statement} we state our generalised cyclic insertion conjecture in an easy, though slightly restricted, form first.
	
	\begin{restatable*}[Generalised cyclic insertion, basic version]{Conj}{conjcyclicinsertioneasy}
		\label{conj:generalcyclicinsertion}
		Let \( B = (0; \ell_1,\ldots,\ell_n) \) be a non-trivial alternating block decomposition of weight \( N \).  Moreover, assume that there is no pair of consecutive lengths \( (\ell_i, \ell_{i+1}) = (1,1) \), including in the case \( (\ell_n, \ell_1) \).  
		\ifprintextra
		Then
		\[
		\sum_{\sigma \in C_n} I_\bl^\mot(\ell_{\sigma(1)},\ldots, \ell_{\sigma(n)}) \eqN I_\bl^\mot(N+2) \, .
		\]
		
		More plainly, using the evaluations of \( I_\bl^\mot(N+2) \) for odd and even \( N \), this conjecture says
		\[
		\sum_{\sigma \in C_n} I_\bl^\mot(\ell_{\sigma(1)},\ldots, \ell_{\sigma(n)}) \eqN \begin{cases}
		(-1)^{N/2} \zeta^\mot(\{2\}^{N/2}) & \text{for weight \( N \) even,} \\
		0 & \text{for weight \( N \) odd.}
		\end{cases}
		\]
		\else
		Then
		\[
		\sum_{\sigma \in C_n} I_\bl(\ell_{\sigma(1)},\ldots, \ell_{\sigma(n)}) \eqN \begin{cases}
		(-1)^{N/2} \zeta(\{2\}^{N/2}) & \text{for weight \( N \) even,} \\
		0 & \text{for weight \( N \) odd.}
		\end{cases}
		\]
		\fi
	\end{restatable*}
	\noindent In \autoref{conj:generalcyclicinsertionfull}, we also present a version with no restriction on the lengths \( (\ell_i, \ell_{i+1}) \).  In \autoref{sec:cyclicinsertion:evidence} we point to \autoref{cor:symins} as evidence for the generalised cyclic insertion conjecture, as well as presenting  numerical evidence for the conjecture.  We then discuss potential ways of tackling this conjecture within the motivic MZV framework.  Unfortunately since the motivic MZV framework cannot yet pin down the rational multiple of \( \zeta^\mot(N) \), a motivic proof by recursion is not possible.  In \autoref{cor:evenweightconditional}, we do manage to offer a proof of the even weight case up to \( \Q \), \emph{conditional} on the odd weight case holding \emph{exactly} at all lower weights.  Then in \autoref{sec:cyclicinsertion:examples} we illustrate carefully how the full version of the conjecture works, and give examples of the identities which it predicts.
	
	In \autoref{sec:123mzv:cyclicinsertion} we introduce the notion of a 123-MZV, and give the statement of the generalised cyclic insertion conjecture for 123-MZV's.  For 123-MZV's, the cyclic insertion identity does not need to be \emph{regularised} before converting back to MZV's so the identities retain more structure; moreover the identities can be generated by a some kind of `local' manipulation of the arguments using the `cyclic operator' \( \cyc \).  In \autoref{sec:123mzv:eg} we present various examples of the cyclic insertion conjecture for 123-MZV's, including the BBBL cyclic insertion conjecture and generalisations of Hoffman's identity.  Then in \autoref{sec:123mzv:differentsym}, we consider some other possible motivic symmetrisations, beyond what is proven with the framework in \autoref{cor:symins}.  These may potentially lead to an even more general framework for motivcally generating identities.  
		
	Finally, in \autoref{sec:furtheridentities} we consider some other unrelated conjectural identities which can understood using the alternating block decomposition, and which generalise a closing conjecture from Section 7.2 of \cite{borwein1998combinatorial}. \medskip
	
	We end this introduction with some informal examples illustrating the prototypical instances of the generalised cyclic insertion conjecture, and how these results follow from the above framework.
	
	\begin{Eg}[BBBL cyclic insertion]\label{sec:intro:eg1}
		The following multiple zeta value
		\[
			\zeta(\{2\}^{b_0}, 1, \{2\}^{b_1}, 3, \ldots, 1, \{2\}^{b_{2n-1}}, 3, \{2\}^{b_{2n}})
		\]
		occurs in the BBBL cyclic insertion conjecture.  It corresponds to the iterated integral
		\[
			\pm I( 0; (10)^{b_0} 1 \mid (10)^{b_1} 10 \mid 0 \cdots 1 \mid (10)^{b_{2n-1}} 10 \mid 0(10)^{b_{2n}}; 1) \, .
		\]
		with the indicated decomposition where alternating blocks are separated by `\( \mid \)', and the sign \( \pm \) is given by the depth of the MZV.  By recording the lengths of these blocks \( \ell_i = 2 b_i + 2 \), we obtain the `block integral' \( \pm I_\bl(2b_0 + 2, 2b_1 + 2, \ldots, 2b_{2n}+2) \). \medskip
		
		By applying \autoref{conj:generalcyclicinsertion} to this correspondence, and using the evaluation \( \zeta(\{2\}^{m}) = \frac{\pi^{\wt}}{(\wt+1)!} \) from \autoref{eq:zeta222}, we recover the BBBL cyclic insertion conjecture
		\[
			\sum_{\sigma \in C_{2n+1}} \zeta(\{2\}^{b_{\sigma(0)}}, 1, \{2\}^{b_{\sigma(1)}}, 3, \ldots, 1, \{2\}^{b_{\sigma(2n-1)}}, 3, \{2\}^{b_{\sigma(2n)}}) \eqN \frac{\pi^\wt}{(\wt+1)!} \, .
		\]
		(In \autoref{eg:bbblcyclic} we work out the sign more carefully, and see the above sign is correct.) \medskip
		
		By applying \autoref{cor:symins}, we obtain
		\[
		\sum_{\sigma \in S_{2n+1}} \zeta(\{2\}^{b_{\sigma(0)}}, 1, \{2\}^{b_{\sigma(1)}}, 3, \ldots, 1, \{2\}^{b_{\sigma(2n-1)}}, 3, \{2\}^{b_{\sigma(2n)}}) \in \pi^\wt \Q \, .
		\]
		This immediately recovers Theorem 3.1 from \cite{charlton2015rational}.
	\end{Eg}

	\begin{Eg}[Hoffman]\label{sec:intro:eg2}
	The following MZV
	\[
		\zeta(3,3,\{2\}^m) 
	\]
	occurs in Hoffman's conjectural identity.  It corresponds to the iterated integral
	\[
		\pm I(0; 10 \mid 010 \mid 0(10)^m ; 1) \, ,
	\] with the indicated decomposition where alternating blocks are separated by `\( \mid \)', and the sign \( \pm \) is given by the depth of the MZV.  By recording the lengths, we obtain the `block integral' \( \pm I_\bl(3,3,2m+2) \).  We also have the correspondences \( -\zeta(3,\{2\}^m,1,2) \leftrightarrow \pm I_\bl(3,2m+2,3) \) and \( \zeta(\{2\}^m,1,2,1,2) \leftrightarrow \pm I_\bl(3,3,2m+2) \), and by \emph{duality} we have \( \zeta(\{2\}^m,1,2,1,2) = \zeta(3,3,\{2\}^m) \). \medskip
	
	By applying \autoref{conj:generalcyclicinsertion} to these correspondences, and using the evaluation \( \zeta(\{2\}^{m}) = \frac{\pi^{\wt}}{(\wt+1)!} \) from \autoref{eq:zeta222}, we recover Hoffman's conjectural identity
	\[
		2 \zeta(3,3,\{2\}^m) - \zeta(3,\{2\}^m,1,2) \eqN -\frac{\pi^\wt}{(\wt+1)!} \, .
	\]
	(In \autoref{eg:hoffman33} we work out the sign more carefully, and see the above sign is correct.)\medskip
	
	By applying \autoref{cor:symins}, obtain a proof of (twice) Hoffman's identity, up to a rational
	\[
		2 \zeta(3,3,\{2\}^m) - \zeta(3,\{2\}^m,1,2) \in \pi^\wt \Q \, .
	\]
	This result is presented in \autoref{thm:hoffman33}, below.
	\end{Eg}

	\subsection{Notation for `equality up to' relations}
	\label{sec:notation}
	
	In this paper, we employ a number of different `equality up to' relations in order to carefully specify the level to which various identities are known to hold.  We use the following
	
	\begin{center}
	\renewcommand{\arraystretch}{1.2}
	\begin{tabular}{c|l}
		Relation & Meaning \\ \hline
		$ \eqN $ & Identity holds numerically to \( \geq 500 \) decimal places \\
		$ \eqQ $ & Identity is proven to hold, up to a rational constant \\
		$ \eqQone $ & The rational constant in \( \eqQ \) is 1, numerically \\ \hline
		$ = $ & Identity is proven to hold exactly
	\end{tabular}
	\end{center}

	\subsection*{Acknowledgements} 
	
	This work was undertaken as part my PhD thesis, and was completed with the support of Durham Doctoral Scholarship funding.  The preparation of this paper occurred while in T\"ubingen, with the support of Teach@T\"ubingen Scholarship funding.
	
	I would like to thank Herbert Gangl for some useful comments on a draft of this paper, which have helped improve the exposition.  I would also like to thank Francis Brown for his comments and feedback on my PhD thesis, in particular on Chapter 2 from which this paper largely derives.  Finally, I would like to thank Eric Panzer for some helpful discussions, and pointing out the result in \autoref{obs:dodd:simplification}.
	
	\section{Background on motivic MZV's and iterated integrals}
	\label{sec:motivicmzvbackground}
	
	\subsection{Properties of motivic MZV's}
		
		In Section 2 of \cite{goncharov2005galois}, Goncharov shows how the classical iterated integrals
		\[
			I(a_0; a_1, \ldots, a_N; a_{N+1})
		\]
		can be lifted to motivic iterated integrals 
		\[
			I^\gmot(a_0; a_1, \ldots, a_N; a_{N+1}) \, , 
		\]
		with new algebraic structure.  This structure comes in the form of a coproduct \( \Delta \), explicitly computed in Theorem 1.2 of \cite[][p. 3]{goncharov2005galois}, making the motivic iterated integrals into a Hopf algebra.  However in Goncharov's setting the motivic MZV \[
			\zeta^\gmot(2) \coloneqq -I^\gmot(0; 1, 0; 1)
		\] vanishes.
		
		In Section 2 of \cite{brown2012mixed}, Brown further lifts the subset of Goncharov's motivic iterated integrals where all \( a_i \in \Set{0, 1} \), in such a way that \( I^\mot(0; 1, 0; 1) \) and the corresponding motivic MZV \( \zeta^\mot(2) \) are now non-zero.  More generally Definition 3.6 of \cite[][p. 8]{brown2012decomposition} defines a motivic MZV as
		\[
		\zeta^\mot(n_1, n_2, \ldots, n_r) \coloneqq (-1)^r I^\mot(0; \underbrace{1, 0, \ldots, 0}_{\text{\( n_1 \) terms}}, \, \underbrace{1, 0, \ldots, 0}_{\text{\( n_2 \) terms}}, \, \ldots \, , \underbrace{1, 0, \ldots, 0}_{\text{\( n_r \) terms}}; 1)
		\, , \]
		in analogy with the Kontsevich integral representation of an MZV, Section 9 in \cite{zagier1994applications}.  Notice here that (motivic) MZV's correspond to iterated integrals with \( a_i \in \Set{0, 1} \), that start with \( a_0 = 0, a_1 = 1 \) and end with \( a_{N} = 0, a_{N+1} = 1 \), namely  having the form \( I^{(\mot)}(0; 1, \ldots, 0; 1) \). \medskip
		
		Brown's motivic MZV's form a graded coalgebra, denoted \( \mathcal{H} \).  The period map
		\begin{equation}\label{eqn:permap}
		\begin{aligned}
		\per \colon \mathcal{H} & \to \R \\
		I^\mot(a_0; a_1, \ldots, a_N; a_{N+1}) &\mapsto I(a_0; a_1, \ldots, a_N; a_{N+1})
		\end{aligned}
		\end{equation}
		defines a ring homomorphism from the graded coalgebra \( \mathcal{H} \) to \( \R \), see Equation 2.11 in \cite[][p. 4]{brown2012mixed} and Equation 3.8 in \cite[][p. 7]{brown2012decomposition}.  This means any identities between motivic MZV's descend to the same identities between ordinary MZV's.\medskip
		
		Theorem 2.4 of \cite[][p. 6]{brown2012mixed} shows that Goncharov's coproduct lifts to a coaction \( \Delta \colon \mathcal{H} \to \mathcal{A} \otimes_\Q\mathcal{H} \) on Brown's motivic MZV's, where \( \mathcal{A} \coloneqq \mathcal{H} / \zeta^\mot(2)\mathcal{H} \) kills \( \zeta^\mot(2) \).  In Section 5 of \cite{brown2012decomposition}, Brown describes an algorithm for decomposing motivic MZV's into a chosen basis using an infinitesimal version of this coaction \( \Delta \colon \mathcal{H} \to \mathcal{A} \otimes_\Q \mathcal{H} \).
		
		The infinitesimal coaction factors through the operators
		\[
		D_r \colon \mathcal{H}_N \to \mathcal{L}_r \otimes_\Q \mathcal{H}_{N-r}
		\, , \]
		where \( \mathcal{L}_r \) is the degree \( r \) component of \( \mathcal{L} \coloneqq \mathcal{A}_{>0} / \mathcal{A}_{>0}\mathcal{A}_{>0} \), the Lie coalgebra of indecomposables, and \( \mathcal{H}_{N} \) is the degree \( N \) component of \( \mathcal{H} \).  The action of \( D_r \) on the motivic iterated integral \( I^\mot(a_0; a_1, \ldots, a_N; a_{N+1}) \) is given explicitly by
		\[
		\sum_{p=0}^{N-r} I^\lmot(a_p; a_{p+1}, \ldots, a_{p+r}; a_{p+r+1}) \otimes I^\mot(a_0; a_1, \ldots, a_p, a_{p+r+1}, \ldots, a_N; a_{N+1})
		\, , \]
		according to Equation 3.4 of \cite[][p. 8]{brown2012mixed}.
		
		The operators \( D_r \) have a pictorial interpretation similar to that of Goncharov's coproduct and the coaction above.  One can view \( D_{r} \) as cutting segments of (interior) length \( r \) out of a semicircular polygon whose vertices are decorated by \( a_0, a_1, \ldots, a_N, a_{N+1} \):
		\begin{center}
			\ifpdf
\begingroup%
  \makeatletter%
  \providecommand\color[2][]{%
    \errmessage{(Inkscape) Color is used for the text in Inkscape, but the package 'color.sty' is not loaded}%
    \renewcommand\color[2][]{}%
  }%
  \providecommand\transparent[1]{%
    \errmessage{(Inkscape) Transparency is used (non-zero) for the text in Inkscape, but the package 'transparent.sty' is not loaded}%
    \renewcommand\transparent[1]{}%
  }%
  \providecommand\rotatebox[2]{#2}%
  \ifx\svgwidth\undefined%
    \setlength{\unitlength}{224bp}%
    \ifx\svgscale\undefined%
      \relax%
    \else%
      \setlength{\unitlength}{\unitlength * \real{\svgscale}}%
    \fi%
  \else%
    \setlength{\unitlength}{\svgwidth}%
  \fi%
  \global\let\svgwidth\undefined%
  \global\let\svgscale\undefined%
  \makeatother%
  \begin{picture}(1,0.5)%
    \put(0,0){\includegraphics[width=\unitlength,page=1]{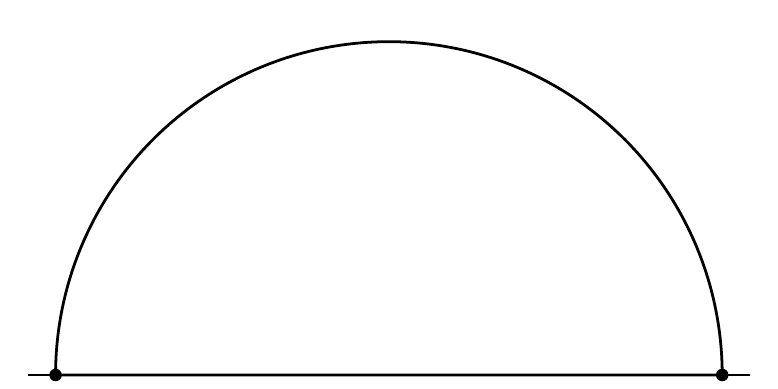}}%
    \put(0.96428574,0.03571433){\color[rgb]{0,0,0}\makebox(0,0)[b]{\smash{$a_N$}}}%
    \put(0.03571429,0.03571433){\color[rgb]{0,0,0}\makebox(0,0)[b]{\smash{$a_0$}}}%
    \put(0.03101021,-0.03617854){\color[rgb]{0,0,0}\makebox(0,0)[b]{\smash{}}}%
    \put(0,0){\includegraphics[width=\unitlength,page=2]{semicircular_dr.pdf}}%
    \put(0.04586007,0.11438756){\color[rgb]{0,0,0}\makebox(0,0)[b]{\smash{$a_1$}}}%
    \put(0.12438499,0.29075758){\color[rgb]{0,0,0}\makebox(0,0)[b]{\smash{$a_{p-1}$}}}%
    \put(0.18933226,0.36288871){\color[rgb]{0,0,0}\makebox(0,0)[b]{\smash{$a_p$}}}%
    \put(0.2678572,0.41994045){\color[rgb]{0,0,0}\makebox(0,0)[b]{\smash{$a_{p+1}$}}}%
    \put(0.64347221,0.45941908){\color[rgb]{0,0,0}\makebox(0,0)[b]{\smash{$a_{p+r}$}}}%
    \put(0.77678738,0.41754635){\color[rgb]{0,0,0}\makebox(0,0)[b]{\smash{$a_{p+r+1}$}}}%
    \put(0.84047309,0.35830366){\color[rgb]{0,0,0}\makebox(0,0)[b]{\smash{$a_{p+r+2}$}}}%
    \put(0.97160691,0.11810024){\color[rgb]{0,0,0}\makebox(0,0)[b]{\smash{$a_{N-1}$}}}%
    \put(0.09216722,0.19943612){\color[rgb]{0,0,0}\rotatebox{66.00000022}{\makebox(0,0)[b]{\smash{$\ldots$}}}}%
    \put(0.43190839,0.46184019){\color[rgb]{0,0,0}\rotatebox{5.99999979}{\makebox(0,0)[lb]{\smash{$\ldots$}}}}%
    \put(0.88661853,0.24107146){\color[rgb]{0,0,0}\rotatebox{-59.99999989}{\makebox(0,0)[lb]{\smash{$\ldots$}}}}%
    \put(0,0){\includegraphics[width=\unitlength,page=3]{semicircular_dr.pdf}}%
  \end{picture}%
\endgroup%

			\else
				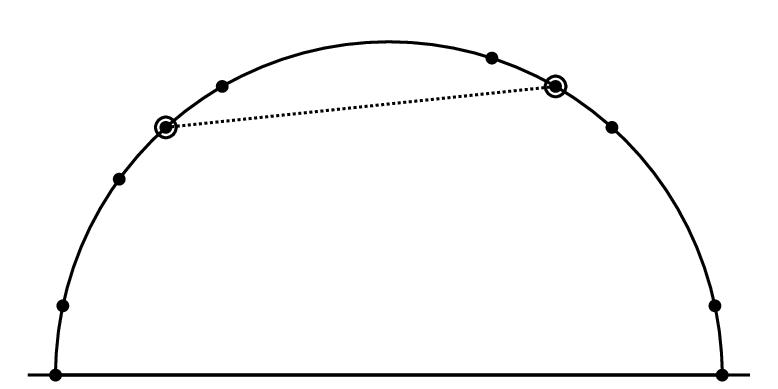
			\fi
		\end{center}
		Notice that the \emph{boundary terms} \( a_p \) and \( a_{p+r+1} \) appear in both the left and right hand factors of \( D_r \), they are part of both the main polygon and the cut-off segment above.
		
		One could also see the operators \( D_r \) as cutting out strings of length \( r+2 \) from the sequence \( (a_0; a_1, \ldots, a_N; a_{N+1}) \).  Following Brown, Definition 4.4 in \cite[][p. 11]{brown2012decomposition}, we make the following definition.
		
		\begin{Def}[Subsequence, and quotient sequence]
			We call the sequence
			\[
			(a_p; a_{p+1}, \ldots, a_{p+r}; a_{p+r+1})
			\]
			appearing in the left factor of \( D_r \) the \emph{subsequence}, and we call the sequence
			\[
			(a_0; a_1, \ldots, a_p, a_{p+r+1}, \ldots, a_N; a_{N+1})
			\]
			appearing in the right factor of \( D_r \) the \emph{quotient sequence} of the original sequence.  Again the boundary terms \( a_p \) and \( a_{p+r+1} \) are part of both the subsequence and the quotient sequence.
		\end{Def}
	
		When decomposing a motivic MZV into a basis, the operator \( D_{2r+1} \) is used to extract the coefficient of \( \zeta^\mot(2r+1) \) as a polynomial in this basis, see Section 5 of \cite{brown2012decomposition}.  The upshot of this comes from Theorem 3.3 of \cite[][p. 9]{brown2012mixed}: 
		\begin{Thm}[Brown]\label{thm:kerDN}
			The kernel of \( D_{<N} \coloneqq \bigoplus_{3 \leq 2r + 1 < N} D_{2r+1} \) is \( \zeta^\mot(N) \Q \) in weight \( N \).
		\end{Thm}
		In other words, if the operators \( D_{2r+1} \), for \( r \) such that \( 3 \leq 2r+1 < N \), all vanish on a given combination of motivic MZV's of weight \( N \), then this combination is a rational multiple of \( \zeta^\mot(N) \).  This will be an important tool in our proof of \autoref{thm:reflclosedid}, the main theorem of this paper. \medskip
		
		Before we continue we should recall a few properties of motivic iterated integrals that will be used in the proof.  See Section 2.4 of \cite{brown2012mixed} for a complete list of properties.  
		
		\begin{Ppty}\label{ppty:propertiesofintegrals}
			The motivic iterated integrals \( I^\mot(a_0; a_1, \ldots, a_N; a_{N+1}) \) satisfy the following list of properties.
			\begin{itemize}
				\item\label{ppty:equalboundaries} (Equal boundaries)  If \( N \geq 1 \), and \( a_0 = a_{N+1} \), then \( I^\mot(a_0; a_1, \ldots, a_N; a_{N+1}) = 0 \),
				\item\label{ppty:reversalofpaths} (Reversal of paths) \( I^\mot(0; a_1, \ldots, a_N; 1) = (-1)^N I^\mot(1; a_N, \ldots, a_1; 0) \), and
				\item\label{ppty:functoriality} (Functoriality) The integral \( I^\mot(a_0; a_1, \ldots, a_N; a_{N+1}) \) is functorial under the map \( t \mapsto 1-t \), so that
				\( I^\mot(a_0; a_1, \ldots, a_N; a_{N+1}) = I^\mot(1-a_0; 1-a_1, \ldots, 1-a_N; 1-a_{N+1}) \).
			\end{itemize}
		
			The reversal of paths property, and the functoriality property, combine to give the \emph{duality} relation on motivic iterated integrals and motivic MZV's, which states
			\begin{itemize}
			\item\label{ppty:duality} (Duality) \( I^\mot(0; a_1, \ldots, a_N; 1) = (-1)^N I^\mot(0; 1-a_N, \ldots, 1-a_1; 1) \).
			\end{itemize}
		\end{Ppty}
	
		Finally, it is useful to introduce the notation of a \emph{trivial subsequence}.
		
		\begin{Def}[Trivial subsequence]\label{def:trivialsubsequence}
			Suppose \( (a_p; a_{p+1}, \ldots, a_{p+r}; a_{p+r+1}) \) is a subsequence, as arising in the computation of some \( D_r \).  If the boundary terms satisfy \( a_p = a_{p+r+1} \), then we call this a \emph{trivial} subsequence.  This is because the corresponding iterated integral vanishes,
			\[
				I^\lmot(a_p; a_{p+1}, \ldots, a_{p+r}; a_{p+r+1}) = 0 \, ,
			\]
			by the \hyperref[ppty:equalboundaries]{equal boundaries} property, so this subsequence does not contribute to \( D_r \).
		\end{Def}
	
		\subsection{Shuffle regularisation of iterated integrals}\label{sec:shufflereg}
		
		The results in this paper, specifically the main theorem (\autoref{thm:reflclosedid}) and its corollaries, and the `full version' of the generalised cyclic insertion conjecture (\autoref{conj:generalcyclicinsertionfull}) apply whether or not the iterated integrals involved are `divergent'.  Recall, an iterated integral \( I^\mot(a_0;a_1, \ldots,a_N;a_{N+1}) \) is called \emph{divergent} if \( a_0 = a_1 \) or \( a_N = a_{N+1} \); these divergent integrals do not correspond directly to multiple zeta values and so must be (shuffle) \emph{regularised} first. \medskip
		
		The procedure for regularising iterated integrals is outlined by \citeauthor{brown2012decomposition}, in Section 2.4 of \cite{brown2012decomposition} for the real-valued MZV's.  A more precise procedure for doing this is given in Section 5.1 of \cite{brown2012decomposition} using relation R2, henceforth called the \emph{divergence relation}.
		
		\begin{Ppty}[{Divergence relation, Relation R2 in \cite[][p. 13]{brown2012decomposition}}]\label{ppty:divergencerelation}
			The iterated integrals \( I^\mot(a_0; a_1, \ldots, a_N; a_{N+1}) \) satisfy the following relation
			\begin{align*}
			& (-1)^k I^\mot(0; \{0\}^k, 1, \{0\}^{n_1-1}, \ldots, 1,\{0\}^{n_r-1}; 1) = \\
			& \sum_{i_1 + \cdots + i_r = k} \tbinom{n_1 -1 + i_1}{i_1} \cdots \binom{n_r -1 + i_r}{i_r} I^\mot(0; 1, \{0\}^{n_1 + i_1 - 1}, \ldots, 1, \{0\}^{n_r + i_r - 1}; 1)
			\, . \end{align*}
		\end{Ppty}
	
		\begin{Rem}
			This relation can be derived by induction, by applying of the shuffle product identity
			\[
			0 = I^\mot(0; 0; 1) I^\mot(0; w; 1) = I^\mot(0; 0 \shuffle w; 1) , \,
			\]
			and the regularisation \( I^\mot(0; 0; 1) \overset{\text{reg}}{=} 0 \) given in property I1 of Theorem 3.5 in \cite{brown2012decomposition}.
		\end{Rem}
	
		\subsubsection*{Regularisation procedure} Given an iterated integral \( I^\mot(a_0; a_1, \ldots, a_N; a_{N+1}) \), the procedure Brown describes to shuffle regularise it and write it as a combination of multiple zeta values is as follows.
		\begin{enumerate}[label=Step \arabic*),leftmargin=2cm]
			\item Use \hyperref[ppty:equalboundaries]{equal boundaries}, and \hyperref[ppty:reversalofpaths]{reversal of paths} to ensure \( a_0 = 0 \), and \( a_{n+1} = 1 \).
			\item Apply the \hyperref[ppty:divergencerelation]{divergence relation} to obtain a linear combination of integrals with \( a_1 = 1 \).
			\item Apply \hyperref[ppty:duality]{duality} to ensure that \( a_{N} = 0 \) in these integrals.
			\item Finally apply the \hyperref[ppty:divergencerelation]{divergence relation} again to ensure \( a_1 = 1 \).
			\item The result can now be written directly in terms of MZV's.
		\end{enumerate}
		
		\begin{Eg}
			To illustrate this procedure, let us regularise \( z = I^\mot(0; 0, 1, 0, 1, 1; 1) \).  We already have \( a_0 = 0 \) and \( a_{N+1} = 1 \), so we can skip step 1.  In step 2, we apply the \hyperref[ppty:divergencerelation]{divergence relation} to get
			\begin{align*}
			z &= I^\mot(0; 0, 1, 0, 1, 1; 1) \\
			&= (-1)^1 \sum_{i_1 + i_2 + i_3 = 1} \binom{1 + i_1}{i_1} \binom{0 + i_2}{i_2} \binom{0 + i_3}{i_3} I^\mot(0; 1, \{0\}^{1 + i_1}, 1, \{0\}^{0+i_2}, 1, \{0\}^{0+i_3}; 1) \\
			&= - 2 I^\mot(0; 1, 0, 0, 1, 1; 1) - I^\mot(0; 1, 0, 1, 0, 1; 1) - I^\mot(0; 1, 0, 1, 1, 0; 1)
			\, . \end{align*}
			In this result, term 3 is already a convergent integral, so in principle we could ignore it for the rest of the procedure.
						
			In step 3 we then apply \hyperref[ppty:duality]{duality}, and obtain
			\[
			z = - 2 I^\mot(0; 0, 0, 1, 1, 0; 1) - I^\mot(0; 0, 1, 0, 1, 0; 1) + I^\mot(0; 1, 0, 0, 1, 0; 1)
			\, . \]
			
			In step 4, we apply the \hyperref[ppty:divergencerelation]{divergence relation} again.  Doing so has no effect on the third term since it is already convergent.  We obtain
			\begin{align*}
			z = {} & -2 \big( I^\mot(0; 1, 0, 0, 1, 0; 1) + 2 I^\mot(0; 1, 0, 1, 0, 0; 1) + 
			3 I^\mot(0, 1, 1, 0, 0, 0; 1) \big) + {} \\
			& \quad {} - \big( -2 I^\mot(0; 1, 0, 0, 1, 0; 1) - 2 I^\mot(0; 1, 0, 1, 0, 0; 1) \big) + {} \\
			& \quad {} + I^\mot(0; 1, 0, 0, 1, 0; 1) \\
			{} = {} & 2 I^\mot(0; 1, 0, 1, 0, 0, 1) + I^\mot(0; 1, 0, 0, 1, 0; 1) + 
			6 I^\mot(0, 1, 1, 0, 0, 0, 1)
			\, . \end{align*}
			
			At last, all of these integrals are convergent, so in step 5 they can be converted directly to MZV's.  We therefore obtain the shuffle regularisation
			\[
 			z \mathrel{\overset{\text{reg}}{=}} 2\zeta^\mot(2, 3) + \zeta^\mot(3, 2) + 6\zeta^\mot(1, 4) \, .
			\]
		\end{Eg}

	\section{Alternating block decomposition}
	\label{sec:alternatingblock}
	
	As explained in the previous section, multiple zeta values correspond to iterated integrals \linebreak \( I(0;w';1) \), where \( w' \) is a word over the alphabet \( \Set{0,1} \) which starts with \( 1 \) and ends with \( 0 \).  For the purposes of computing the derivations \( D_{2r+1} \), the upper and lower limits of the integral are equally as important as the rest of the word \( w' \).  We therefore want to prioritise the word \( w = 0w'1 \), instead.  We are going to decompose this \( w \) into so-called \emph{alternating blocks}.
	
	\begin{Def}[Alternating block]
		Let \( W_0 \) denote the infinite string \( W_0 \coloneqq 010101\ldots \), consisting of an alternating sequence of 0's and 1's, beginning with a 0.  Let \( W_1 \) denote the infinite string \( W_1 \coloneqq 101010\ldots \), consisting of an alternating sequence of 1's and 0's, beginning with a 1.
		
		We write \( W_i^\ell \) to denote the string obtained by taking the first \( \ell \) letters of \( W_i \).  Such a finite string is called an \emph{alternating block}.
	\end{Def}
	
	It is clear that every word over the alphabet \( \Set{0,1} \) can be written as a concatenation of various alternating blocks \( W_{\epsilon_i}^{\ell_i} \), simply because \( W_0^1 = 0 \) and \( W_1^1 = 1 \).  We are interested in decompositions containing the minimal number of blocks \( \size{\Set{W_{\epsilon_i}^{\ell_i}}} \), which will be obtained by `cutting' at a repeated letter.
	
	\begin{Lem}\label{lem:altblockdecomp}
		Every word \( w \) over \( \Set{0,1} \) can be decomposed as a concatenation \( W_{\epsilon_1}^{\ell_1} \cdots W_{\epsilon_n}^{\ell_n} \), such that the last letter of \( W_{\epsilon_i}^{\ell_i} \) agrees with the first letter of \( W_{\epsilon_{i+1}}^{\ell_{i+1}} \).  Moreover this representation is unique, and involves the minimal number \( n = \size{\Set{W_{\epsilon_i}^{\ell_i}}} \) of blocks.
		
		\begin{proof}  We deal separately with the existence, uniqueness and minimality claims.
			
			\proofstep{Existence:} Such a decomposition exists by induction.  For words of length \( 1 \), we have
			\[
				w = \epsilon \text{ decomposes as } W_\epsilon^1 \, .
			\]
			
			Now suppose all words of length \( < L \) can be so decomposed, and let \( w \) be a word of length \( L \).  If \( w \) does not contain the substring \( 00 \), or the substring \( 11 \), \( w \) is already alternating, so \( w = W_{\epsilon_1}^L \), where \( \epsilon_1 \) is the first digit of \( w \).
			
			Otherwise find the first occurrence of \( 00 \) or \( 11 \) in \( w \), and split \( w = w_1 w_2 \) at this point.  Both \( w_1 \) and \( w_2 \) have length \( < L \), so can be decomposed as required. Say \( w_1 = W_{\epsilon_1}^{\ell_1}\cdots W_{\epsilon_n}^{\ell_n} \), and \( w_2 = W_{\epsilon_{n+1}}^{\ell_{n+1}} \cdots W_{\epsilon_{n+m}}^{\ell_{n+m}} \).  By construction the last digit of \( W_{\epsilon_n}^{\ell_n} \) agrees with the first digit of \( W_{\epsilon_{n+1}}^{\ell_{n+1}} \), so concatenating these decompositions gives a decomposition of \( w \) of the required form.
			
			\proofstep{Uniqueness:} Suppose \( B_1 = W_{\epsilon_1}^{\ell_1} \cdots W_{\epsilon_n}^{\ell_n} \) and \( B_2 = W_{\delta_1}^{k_1} \cdots W_{\delta_m}^{k_m} \) are two ostensibly different decompositions for \( w \).
			
			We can remove any leading terms \( B_0 = W_{\epsilon_1}^{\ell_1} \cdots W_{\epsilon_m}^{\ell_m} \) which happen to agree.  If this produces two copies of the empty word, then we are done, so assume not.  If only one copy of the empty word results, say \( B_1 = B_0 \), and \( B_2 = B_0 B_2' \), then \( B_1 \) and \( B_2 \) have different lengths, so cannot be decompositions of the same word \( w \).
			
			Hence after removing any common leading terms, we obtain two non-empty words \( B_1 = B_0 B_1' \) and \( B_2 = B_0 B_2' \).  We can therefore assume that the first difference between \( B_1 \) and \( B_2 \) occurs in the first term.
			
			Since \( B_1 \) and \( B_2 \) describe the same word \( w \), we necessarily have \( \epsilon_1 = \delta_1 \).  Since the first terms differ, we must have \( \ell_1 \neq k_1 \), say \( \ell_1 < k_1 \)
			
			We can now determine the letter at position \( \ell_1 + 1 \).  Using \( B_1 \) it is \( \epsilon_1 + (\ell_1 - 1) \pmod*{2} \), by using \( B_2 \) it is \( \epsilon_1 + (\ell_1) \pmod*{2} \).
			\begin{align*}
				B_1 &\rightsquigarrow (\underbrace{0101\cdots 01}_{\text{\(\ell_1\) symbols}}) (\underbrace{1}_{\mathrlap{\hspace{-10pt}\text{position \( \ell_1 + 1 \)}}}0 \cdots 10) \cdots \\
				B_2 &\rightsquigarrow (\underbrace{0101\cdots 01}_{\text{\( \ell_1 \) symbols}} \underbrace{0}_{\mathrlap{\hspace{-10pt}\text{position \( \ell_1 + 1 \)}}} \cdots) \cdots \, .
			\end{align*}
			These are not equal, contradicting the fact that \( B_1 \) and \( B_2 \) describe the same word.
			
			Thus is it not possible to have two different decompositions for the word \( w \).
			
			\proofstep{Minimality:} First observe \( \size{(\text{00's or 11's in \( w \)})} \leq \size{\text{blocks}} - 1 \). Each block is alternating, so the only place \( 00 \) or \( 11 \) can occur is where two blocks join, and there are \( \size{\text{blocks}}-1 \) such places.  Alternatively, we can read this as \( \size{\text{blocks}} \geq \size{(\text{00's or 11s in \( w \)})} + 1 \), giving a lower bound for the number of blocks.
			
			With the construction above we do have 00 or 11 at each point where two block join, so \( \size{(\text{00's or 11's in \( w \)})} = \size{\text{blocks}} - 1 \), giving equality.  Hence the lower bound for the number of blocks in the decomposition is obtained.
		\end{proof}
	\end{Lem}
	
	From this we define the alternating block decomposition of a word \( w \).
	
	\begin{Def}[Alternating block decomposition]
		Let \( w \) be a word over \( \Set{0,1} \), and let \( W_{\epsilon_1}^{\ell_1} \cdots W_{\epsilon_n}^{\ell_n} \) be the decomposition from \autoref{lem:altblockdecomp}.  The \emph{(alternating) block decomposition} of \( w \) is defined to be
		\[
			\block(w) \coloneqq (\epsilon_1; \ell_1,\ldots,\ell_n) \, ,
		\]
		obtained by recording the starting letter \( \epsilon_1 \) of \( w \), and the lengths \( \ell_i \) of each of the alternating blocks.
	\end{Def}
	
	\begin{Rem}
		Notice that only \( \epsilon_1 \) out of the \( \epsilon_i \) is required in this description.  We can calculate \( \epsilon_{i+1} \) from \( \epsilon_i \) by knowing that the first digit of \( W_{\epsilon_{i+1}}^{\ell_{i+1}} \) is equal to the last digit of \( W_{\epsilon_i}^{\ell_i} \).  We find that \( \epsilon_{i+1} = \epsilon_i + (\ell_i - 1) \pmod*{2} \).
		
		We can therefore recover \( w \) from \( B \coloneqq \block(w) \),  by computing the other \( \epsilon_i \), and finding \( \word(B) \coloneqq W_{\epsilon_1}^{\ell_1} \cdots W_{\epsilon_n}^{\ell_n}. \)
	\end{Rem}
	
	\begin{Eg}\label{eg:blockdecomposition}
		We have the following block decompositions.  We indicate where the blocks occur in \( w \) by writing `\( \mid \)'.
		\begin{alignat*}{4}
			& w_1 = 01010 \mid 01 \mid 1 \mid 1010101 & \text{ has } & \block(w_1) = (0; 5, 2, 1, 7) \\
			& w_2 = 1 \mid 101 \mid 1010 \mid 0 \mid 01010 & \text{ has } & \block(w_2) = (1; 1, 3, 4, 1, 5) \, .
		\end{alignat*}
	\end{Eg}
	
	\begin{Not}
		Given a block decomposition \( B = (\epsilon_1; \ell_1,\ldots,\ell_n) \), we define the associated (motivic) \emph{block integral} as
		\[
			I^{(\mot)}_\bl(B) = I^{(\mot)}_\bl(\epsilon_1; \ell_1,\ldots,\ell_n) \coloneqq I^{(\mot)}(\word(B)) \, .
		\]
		
		In the case were \( \epsilon_1 = 0 \), we may write \( B = (\ell_1,\ldots,\ell_n) \), and \( I_\bl(B) = I_\bl(\ell_1,\ldots,\ell_n) \) for ease of notation.
	\end{Not}

	\begin{Eg}
		Using the block decompositions from \autoref{eg:blockdecomposition}, we can write the corresponding iterated integrals as
		\begin{alignat*}{4}
			I^\mot(w_1) &= I^\mot(0; 1, 0, 1, 0, 0, 1, 1, 1, 0, 1, 0, 1, 0; 1) &&= I_\bl^\mot(5, 2, 1, 7) \\
			I^\mot(w_2) &= I^\mot(1; 1, 0, 1, 1, 0, 1, 0, 0, 0, 1, 0, 1; 0) &&= I_\bl^\mot(1; 1, 3, 4, 1, 5) \, .
		\end{alignat*}
	\end{Eg}
	
	\begin{Not}
		Given a block decomposition \( B = (\epsilon_1; \ell_1, \ldots, \ell_n) \) we introduce the following notation (in terms of \( B \) \emph{only}): the length of the \( i \)-th block is \( B_i^\len \); the starting digit of the \( i \)-th block is \( B_i^\st \); and the ending digit of the \( i \)-th block is \( B_i^\en \).
	\end{Not}
	
	Now we collect some simple facts about these block decompositions, and block integrals.
	
	\begin{Lem}
		The block integral \( I_\bl^\mot(\epsilon_1; \ell_1,\ldots,\ell_n) \) has weight \( -2 + \sum_i \ell_i \).
		
		\begin{proof}
			The word \( w = \word(\epsilon_1; \ell_1,\ldots,\ell_n) \) has length \( \sum_i \ell_i \) because the \( i \)-th block has length \( \ell_i \).  This includes the lower bound and upper bound of the integral \( I^\mot(w) \), so subtract 2.
		\end{proof}
	\end{Lem}
	
	We use this lemma to give meaning to the notion of the \emph{weight} of a block decomposition.

	\begin{Lem}\label{lem:weightandblockstrivial}
		Let \( I = I_\bl^\mot(B) \) be a block integral with weight \( N \), and with \( n \) blocks.  The upper and lower bounds of \( I \) are equal, meaning \( I = 0 \), if and only if \( N = n \pmod*{2} \).
		
		\begin{proof}
			We can compute the \( i \)-th digit (counting from 1) of the alternating string \( W_{\epsilon_i} \) to be \( \epsilon_1 + (i-1) \pmod*{2} \).
			
			Each new block in a block decomposition adds a delay of 1: what should be a 0 is now a 1, and vice versa.  So the last digit of \( \word(B) \) is \( \epsilon_1 + (\sum_i \ell_i - 1) - (n - 1) = \epsilon_1 + t - n \pmod*{2} \).  This is equal to \( \epsilon_1 \) if and only if \( t - n = 0 \pmod*{2} \), giving the required result.
		\end{proof}
	\end{Lem}

	\begin{Rem}
		It is perhaps interesting to compare the structure of this with Tsumura's depth-parity	theorem \cite{tsumura2004combinatorial}. Both results have that form that an object simplifies (to zero in this case, or to lower depth in Tsumura's case) if some parity condition holds (equal parity in this case, and opposite parity	in Tsumura's case).
		
		Of course Tsumara's result is highly non-trivial fact, whereas this result is just a trivial observation that the bounds of integration are equal, 
	\end{Rem}
	
	If, for some block decomposition \( B \), it is the case that \( N = n \pmod*{2} \), we will say the block decomposition is \emph{trivial}.
	
	\begin{Lem}\label{lem:divergent}
		Let \( I = I_\bl^\mot(B) \), where \( B = (\epsilon_1; \ell_1,\ldots,\ell_n) \) is a non-trivial block decomposition, with \( n \) blocks and weight \( N > 2 \).  The integral \( I \) is divergent if and only if \( \ell_1 = 1 \) or \( \ell_n = 1 \).
		
		\begin{proof}
			If \( \ell_1 = 1 \), then the integral \( I \) begins \( \epsilon_1\epsilon_1 \cdots \), so is divergent.  If \( \ell_n = 1 \), then the integral \( I \) ends \( \cdots \epsilon_n\epsilon_n \), so is also divergent.
			
			Otherwise \( \ell_1 > 1 \), and \( \ell_n > 1 \) meaning the integral starts and ends \( \epsilon_1 (1-\epsilon_1) \cdots \epsilon_n (1 - \epsilon_n) \), so is convergent.
		\end{proof}
	\end{Lem}
	
	\section{Reflection operators, and reflectively closed sets}
	\label{sec:reflectionoperators}
	
	Ultimately, we wish to compute the operators \( D_{2r+1}\) on certain combinations of motivic iterated integrals arising from block decompositions, and see \( D_{2r+1} \) vanishes.  We do this via reflection operators defined on subsequences and on block decompositions.
	
	\begin{Def}[Reflection \( \refl_{j,k} \)]
		Let \( B = (\epsilon_1; \ell_1,\ldots,\ell_n) \) be a (non-trivial) block decomposition with \( n \) blocks.  For each \( 1 \leq j \leq k \leq n \), we define the \emph{reflection operator} \( \refl_{j,k} \) as follows.
		
		We set \( \refl_{j,k} B \coloneqq (\epsilon_1'; \ell_1', \ldots, \ell_n') \), where \( \epsilon_1' = \epsilon_1 \), and
		\[
			\ell_i' \coloneqq \begin{cases}
				\ell_i & \text{for \( i < j \), or \( i > k \), and } \\
				\ell_{k+j-i} & \text{for \( j \leq i \leq k \),} 
			\end{cases}
		\]
		that is to say, we reverse the lengths of the blocks from position \( j \) to position \( k \), inclusive.  The reflection operator \( \refl_{j,k} \) is extended to \emph{words} by using \( \refl_{j,k} w = \word(\refl_{j,k} \block(w)) \).
	\end{Def}
	
	The following results are immediate consequences of this definition.
	
	\begin{Lem}
		Let \( B \) be a block decomposition with \( n \) blocks.  Then for any \( 1 \leq j \leq k \leq n \), the reflection operator \( \refl_{j,k} \) preserves the weight of \( B \), and the number of blocks in \( B \).
	\end{Lem}
	
	\begin{Lem}
		Let \( X \) be the set of all block decompositions with some given number \( n \) of blocks.  Then for any \( 1 \leq j \leq k \leq n \), the reflection operator \( \refl_{j,k} \) is an involution on \( X \), so that \( \refl_{j,k}^2 B = B \) for any \( B \in X \).
	\end{Lem}

	\begin{Eg}
		The word \( w_1 = 01010 \mid 01 \mid 1 \mid 1010101 \) from \autoref{eg:blockdecomposition}, has block decomposition \( B_1 = \block(w_1) = (5, 2, 1, 7) \).  We compute that
		\[
			\refl_{1,3} B_1 = \refl_{1,3} (5, 2, 1, 7) = (1, 2, 5, 7) \, .
		\]
		In terms of words, we can write
		\begin{align*}
			\refl_{1,3} w_1 & {} = \refl_{1,3} 01010 \mid 01 \mid 1 \mid 1010101 \\
				& {} = 0 \mid 01 \mid 10101 \mid 1010101 \, .
		\end{align*}
	\end{Eg}
	
	\begin{Rem}\label{rem:dualityr1n}
		The reflection operator \( \refl_{1,n} \), on a non-trivial block decomposition with weight \( N \) and \( n \) blocks, is closely related to the duality of (motivic) iterated integrals.  The operator \( \refl_{1,n} \) reverses all of the block lengths.  Since the word starts 0 and ends 1, this has the effect of reversing the word and interchanging \( 0 \leftrightarrow 1 \).  So we can express duality as \( I_\bl^{(\mot)}(B) = (-1)^N I_\bl^{(\mot)}(\refl_{1,n} B) \), or in terms of words
		\[
			I^{(\mot)}(w) = (-1)^N I^{(\mot)}(\refl_{1,n} w) \, .
		\]
	\end{Rem}
	
	\begin{Rem}
		Although the \( \refl_{j,k} \) are defined on block decompositions and words, they are not well-defined on the associated block integrals.  This is because the block integrals satisfy relations that, amongst other things, do \emph{not} preserve the number of blocks.  For example, the Zagier-Broadhurst identity
		\[
			\zeta(\{1,3\}^n) = \frac{1}{2n+1} \zeta(\{2\}^{2n})
		\]
		corresponds to the following identity on iterated integrals
		\[
			I(0; \{1,1,0,0\}^n; 1) = \frac{1}{2n+1} I(0; \{1,0\}^{2n}; 1) \, .
		\]
		After finding the block decomposition of the words \( 0(1100)^n1 \) and \( 0(10)^{2n}1 \) inside the iterated integrals, we can write the identity in terms of block integrals, as follows
		\[
			I_\bl(\{2\}^{2n+1}) = \frac{1}{2n+1} I_\bl(4n + 2) \, .
		\]
		The integral on the left consists of \( 2n + 1 \) blocks, whereas the integral on the right contains only 1 block.  We see, therefore, that the operator \( \refl_{1,2} \) is not even defined on the 1 block integral, on the right hand side.
		
		Any application of \( \refl_{j,k} \) directly to iterated integrals would be via abuse of notation, and would require acting on the particular \emph{argument string} inside the integral.
	\end{Rem}
	
	We can now define the main objects which will be used to create identities on iterated integrals and MZV's.
	
	\begin{Def}[Reflectively closed set]
		Let \( S \) be a set of block decompositions with fixed weight \( N \), and fixed number of blocks \( n \).  We say that \( S \) is \emph{reflectively closed} if for every \( B \in S \), the result of \( \refl_{j,k} B \) is already in \( S \), for \( 1 \leq j \leq k \leq n \).
	\end{Def}
	
	\begin{Def}[Reflective closure]
		Let \( S \) be a set of block decompositions as in the previous.  We define the \emph{reflective closure} of \( S \), written \( \langle S \rangle_{\refl} \), to be the smallest reflectively closed set containing \( S \).
	\end{Def}
	
	The set of all block decompositions with some fixed weight \( N \), and some fixed number of blocks \( n \) is finite, since there are only finitely many \emph{words} over \( \Set{0,1} \) with length \( N+2 \).  It is easy to see that the intersection of any two reflectively closed sets is also reflectively closed.  So the reflective closure of any set \( S \) certainly exists.  More useful is the form this reflective closure takes.
	
	\begin{Prop}
		Let \( B = (\epsilon_1; \ell_1,\ldots,\ell_n) \) be a block decomposition.  Then the reflective closure of \( \Set{B} \) is given by
		\[
			\langle B \rangle_{\refl} = \Set{ (\epsilon_1; \ell_{\sigma(1)}, \ldots, \ell_{\sigma(n)}) | \sigma \in S_n } \, .
		\]
		
		\begin{proof}
			Certainly \( \Set{ (\epsilon_1; \ell_{\sigma(1)}, \ldots, \ell_{\sigma(n)}) | \sigma \in S_n} \) is reflectively closed, because the reflection operators merely permute the lengths.  Moreover, the operator \( \refl_{i,i+1} \) induces the transposition \( (i \, i+1) \) on the \( \ell_i \), so we necessarily generate every permutation.
		\end{proof}
	\end{Prop}
	
	We consider now the subsequences of a word \( w \), as would appear in the computation of \( D_{2r+1} I(w) \).  We will encode these subsequences using the block decomposition, and define reflection operations on them using the \( \refl_{j,k} \).
	
	\begin{Def}[Encoding of a subsequence]\label{def:subsequence}
		Let \( w \) be a word over \( \Set{0,1} \), and let \( P \) be a subsequence of \( w \) of length \( \geq 2 \), in the sense of the derivations \( D_r \).  We encode the subsequence \( P \) on \( w \) by giving the following data:
		\begin{itemize}
			\item The block encoding \( B = \block(w) \), of the word \( w \) upon which \( P \) is defined,
			\item the block \( s \) in which \( P \) starts,
			\item the block \( t \) in which \( P \) finishes,
			\item the number of letters \( \alpha \) before \( P \), in the block \( s \), and
			\item the number of letters \( \beta \) after \( P \), in the block \( t \).
		\end{itemize}
		
		We package these into the following tuple, and identify it with the subsequence \( P \) to write
		\[
			P = (B; s, t; \alpha, \beta) \, .
		\]
	\end{Def}

	\begin{Eg}\label{eg:subsequenceencoding}
		Consider the indicated subsequence \( P \) on the word \( w_2 = 1 \mid 101 \mid 1010 \mid 0 \mid 01010 \) from \autoref{eg:blockdecomposition}
		\[
			w_2 = 1 \mid \overbrace{1}^{\mathclap{\text{1 letter}}}\underbracket{01 \mid 1010 \mid 0 \mid 01}_{\text{subsequence \( P \)}}\overbrace{010}^{\text{3 letters}}
		\]
		It starts in block 2 of \( w_2 \), and ends in block 5.  It has 1 letter before it in block 2, and it has 3 letters after it in block 5.  It is encoded as
		\[
			P = (\block(w_2); 2, 5; 1, 3) = \big( (1; 1,3,4,1,5); 2, 5; 1, 3 \big) \, .
		\]
	\end{Eg}
	
	From these data we can calculate the length of the subsequence \( P \) to be \( \sum_{i=s}^t B_i^\len - \alpha - \beta \).  We have the following natural restrictions which completely characterise valid subsequences.
	
	\begin{Lem}
		An encoding \( (B; s, t; \alpha, \beta) \) of a subsequence is valid (corresponds to an actual subsequence of length \( \geq 2 \)) if and only if the following conditions hold
		\begin{enumerate}
			\item \( 1 \leq s \leq t \leq n \), where \( n \) is the number of blocks in \( B \),
			\item \( 0 \leq \alpha < B_s^\len \),
			\item \( 0 \leq \beta < B_t^\len \), and
			\item if \( s = t \), then \( \alpha + \beta + 2 \leq B_s^\len \).
		\end{enumerate}
		
		\begin{proof}
			These conditions are necessary for the following reasons.  Item i) corresponds to the fact that a subsequence starts before it finishes, and lies within the word \( w = \word(B) \).  Item ii), the subsequence may start at the first letter of the block (0 letters before it), or the last letter of the block (\( B_s^\len - 1 \) letters before it).  Item iii), similarly the subsequence may end at the last letter of the block, or the first letter of the block.  Item iv), the subsequence lies entirely within the one block, starts before it finishes, and has length \( \geq 2 \).
			
			Conversely, given these conditions, one can uniquely mark out a subsequence on \( w = \word(B) \).  Find blocks \( s \) and \( t \).  Count \( \alpha \) letters from the start of block \( s \) to find the start of the subsequence.  Count \( \beta \) letters from the end of block \( t \) to find the end of the subsequence.  If \( s \neq t \), the subsequence contains points from two different blocks, so has length \( \geq 2 \).  Otherwise condition iv) ensures length \( \geq 2 \).
		\end{proof}
	\end{Lem}
	
	We can now define a reflection operator on subsequences, using this encoding.
	
	\begin{Def}[Reflection of a subsequence]
		Let \( P = (B; s, t; \alpha, \beta) \) be a subsequence on some word \( w = \word(B) \).  The reflection operator \( \refl \) is defined on \( P \) by
		\[
			\refl P \coloneqq (\refl_{s,t} B; s, t; \beta, \alpha) \, .
		\]
	\end{Def}
	
	\begin{Rem}
		Notice that \( s \) and \( t \) remain in the same order, whereas \( \alpha \) and \( \beta \) are reversed by \( \refl \).  There are two reasons for this.  Firstly to be a valid subsequence, we must have \( s \leq t \), so we \emph{cannot} switch the order of \( s \) and \( t \).
		
		Secondly, we want to \( \refl \) to capture the idea of \emph{reflecting} the blocks which contain the subsequence \( P \) and carrying the subsequence \( P \) along in the process; this will not change which block the subsequence starts in, or which block is ends in.
	\end{Rem}

	\begin{Eg}
		Applying the subsequence reflection operator \( \refl \) to the subsequence \( P = \big( (1; 1,3,4,1,5); 2, 5; 1, 3 \big) \) from \autoref{eg:subsequenceencoding} produces the following.
		\begin{align*}
			\refl \big( (1; 1,3,4,1,5); 2, 5; 1, 3 \big) & {} = \big( \refl_{2,5} (1; 1,3,4,1,5); 2,5; 3,1) \\
				&= \big( (1; 1,5,1,4,3); 2, 5; 3, 1 \big) \, .
		\end{align*}
		We can see the `reflection' that \( \refl \) produces by realising \( P \) and \( \refl P \) as subsequences on words.  We have
		\begin{align*}
			P &\rightsquigarrow 1 \mid \overbrace{\underbrace{1}_{\mathllap{\text{1 letter \hspace{-1em}}}}\underbracket{01 \mid 1010 \mid 0 \mid 01}_{\text{subsequence \( P \)}}\underbrace{010}_{\mathrlap{\text{\hspace{-1em} 3 letters}}}}^{\text{reflect}} \\
			& \hspace{8em} \updownarrow \\
			\refl P &\rightsquigarrow 1 \mid \overbrace{\underbrace{101}_{\mathllap{\text{3 letters \hspace{-1em}}}} \underbracket{01 \mid 1 \mid 1010 \mid 01}_{\text{subsequence \( \refl P \)}} \underbrace{0}_{\mathrlap{\text{\hspace{-1em} 1 letter}}}}^{\text{reflect}} \, .
		\end{align*}
		Notice that, as words, the subsequence \( \refl P \) is \emph{dual} to the subsequence \( P \).  We show in \autoref{lem:revordual} that \( \refl P \) is always the reverse of \( P \), or the dual or \( P \).
	\end{Eg}
	
	The following lemmas are immediate consequences of these definitions, and of the conditions characterising a valid subsequence.
	
	\begin{Lem}\label{lem:reflofsub}
		Applying the reflection operator \( \refl \) to the subsequence \( P = (B; s, t; \alpha, \beta) \) produces a valid subsequence on \( \word(\refl_{s,t} B) \).
	\end{Lem}
	
	\begin{Lem}
		The operator \( \refl \) preserves the length of a subsequence.
	\end{Lem}
	
	\begin{Lem}\label{lem:reflinvolution}
		The operator \( \refl \) is an involution on the set of all subsequences on some reflectively closed set containing block decompositions with fixed weight \( N \), and fixed number of blocks \( n \).
	\end{Lem}
	
	\section{Identities from reflectively closed sets}
	\label{sec:identitiesfromreflection}
	
	We aim now to prove the following theorem.
	
	\thmreflectivelyclosed*
	
	\begin{proof}[Proof idea]\renewcommand{\qedsymbol}{}
		The main idea of the proof is that the subsequence reflection operator \( \refl \) sets up a pairwise cancellation between all subsequences on integrals from \( S \).  This means \( D_{<N} \) vanishes, and the result follows by Brown's characterisation of \( \ker D_{<N} \) from \autoref{thm:kerDN}.
	\end{proof}

	We will break the proof into a number of steps, and present it via a sequence of simple lemmas.  To start we need a few technical results.
	
	The following results allows us to discard some subsequences in the computation of \( D_{<N} \).
	
	\begin{Lem}\label{lem:blockrefl_fix}
		Suppose that \( B \) is the block decomposition of some iterated integral, and further suppose that \( \refl_{j,k} B = B \), for some \( j, k \).  Then, in particular, \( B_{k-i}^\len = B_{j+i}^\len \), for \( 0 \leq i \leq k-j \).  Moreover, if \( k-j + 1 \) is even, or if \( k-j + 1 \) is odd and \( B_{j + (k-j)/2}^\len \) is odd, then \( B_j^\st = B_k^\en \).
		
		\begin{proof}
			By the definition of \( \refl_{j,k} \) on \( B \), it is clear that \( B_{j+i}^\len = B_{k-i}^\len \) since \( B_{j+i}^\len = (\refl_{jk} B)_{j+i}^\len = B_{k + j - (j+i)}^\len = B_{k-i}^\len \).
			
			By removing blocks \( < j \), and removing blocks \( > k \), we can assume that the computation is of \( \refl_{1,n} B \), with \( B \) having \( n \) blocks.  The case where \( k-j + 1 \) is even corresponds to \( n \) even, and the case \( k-j + 1 \) is odd corresponds to \( n \) odd.
			
			For \( n = 2m \) even, we have that \( B_{m}^\en = B_{m+1}^\st \) by the definition of a block decomposition.  Since \( B_m^\len = B_{m+1}^\len \), we see that \( B_m^\st = B_{m+1}^\en \).  Continue this outwards until we get \( B_1^\st = B_n^\en \).

			For \( n = 2m + 1 \) odd, we have \( (\refl_{1,n} B)_{m+1}^\st = B_{m+1}^\st \) by assumption.  We have that \( (j + (k-j)/2) \) corresponds to \( 1 + (2m + 1 - 1)/2 = m + 1 \), so that \( B_{m+1}^\len \) is odd.  This means that \( B_{m+1}^\en = B_{m+1}^\st + (B_{m+1}^\len - 1) = B_{m+1}^\st \pmod*{2} \).  Then use the argument above to work outwards to get \( B_1^\st = B_n^\en \). 
		\end{proof}
	\end{Lem}
	
	\begin{Lem}\label{lem:subseqfix_trivial}
		Suppose that the subsequence \( P \) is a fixed point of the reflection operator \( \refl \).  Further, suppose that \( P \) has odd length.  Then \( P \) is trivial.  
		
		\begin{proof}
			If \( P = (B; s, t; \alpha, \beta) \) is a fixed point, then we must have \( \refl_{s,t} B = B \), and \( \alpha = \beta \) by the definition of \( \refl \).
			
			Firstly we show that it is not possible for \( P \) to have odd length, be a fixed point, and have \( t - s + 1 \) even.  For if this were the case, by \autoref{lem:blockrefl_fix} we necessarily have \( B_{s+i}^\len = B_{t - i}^\len \), and no `middle block'.  This means \( P \) has length \( \sum_{i=s}^t B_{i}^\len - \alpha - \beta = 2\sum_{i=s}^{s + (t-s)/2} B_i^\len - 2\alpha = 0 \pmod*{2} \).
			
			Therefore we are in the case where \( t - s + 1 \) is odd.  Here we claim that we must have \( B_{s + (t-s)/2}^\len \) odd.  Otherwise as before, \( P \) would have length \( \sum_{i=s}^t B_{i}^\len - \alpha - \beta = 2\sum_{i=s}^{s + (t-s)/2} B_i^\len + B_{s + (t-s)/2} - 2\alpha = 0 \pmod*{2} \).
			
			Now we can apply \autoref{lem:blockrefl_fix} to conclude that \( B_s^\st = B_t^\en \), for the subsequence \( P \).  Therefore the first digit of \( P \) is \( B_s^\st + \alpha \), whilst the last digit of \( P \) is \( B_t^\en - \beta = B_s^\st - \alpha = B_s^\st + \alpha \pmod*{2} \).  Thus the subsequence is trivial.
		\end{proof}
	\end{Lem}
	
	Applying the reflection operator to subsequences produces the reverse, or the dual word.
	
	\begin{Lem}\label{lem:revordual}
		Let \( P = (B; s, t; \alpha, \beta) \) be a subsequence.  Then, as words, the subsequence \( \refl P \) is either the reverse of \( P \), or the dual of \( P \), i.e. the reverse with \( 0 \leftrightarrow 1 \).
		
		\begin{proof}
			By removing the blocks \( < s \), and the blocks \( > t \), we may assume \( s = 1 \) and \( t = n \), where \( n \) is the number of blocks in \( B \).
			
			In the case where \( B_1^\st = B_n^\en \), we will show that the subsequence \( \refl P \) is the reverse of the subsequence \( P \).  Let \( B = (\epsilon_1; \; \ell_1, \ldots, \ell_n) \).  The first digit of \( P \) is then \( \epsilon_1 + \alpha \).  And \( \refl_{1,n} B = (\epsilon_1; \; \ell_n, \ldots, \ell_1) \), so the last digit of \( \refl P \) is \( (\refl_{1,n}B)_n^\en - \alpha \).  As in the proof of \autoref{lem:weightandblockstrivial}, we compute that \( (\refl_{1,n}B)_n^\en = \epsilon_1 + \sum_{i=1}^n (\ell_i - 1) = B_n^\en \pmod*{2} \) and by assumption this is \( = B_1^\st \pmod*{2} \).  So the last digit of \( \refl P \) is \( B_1^\st - \alpha \), which equals the first.
			
			We can repeat this one letter at a time to see that the subsequence \( \refl P \) is exactly the reverse of \( P \).
			
			In the case where \( B_1^\st \neq B_n^\en \), the \( \refl P \) is the dual of the subsequence \( P \).  Observe that in this case we have that the last digit is \( 1 - B_1^\st - \alpha \), so at every point we have the extra step of taking \( 1 - B_i^\en \).  Not only is the subsequence reversed, but we also interchange \( 0 \leftrightarrow 1 \), giving the dual overall.
		\end{proof}
	\end{Lem}
	
	In what follows, let \( S \) be a reflectively closed subset of block decompositions with some fixed weight \( N \), and some fixed number of blocks \( n \).  Let \( T \) be the set of all odd length subsequences on the block decompositions in \( S \).
	
	\begin{Lem}
		The reflection operator \( \refl \) defines a map from \( T \to T \).
		
		\begin{proof}
			Let \( P \) be a subsequence in \( T \); then \( P =  (B; s, t; \alpha, \beta) \) for \( B \) a block decomposition in \( S \), and some \( s, t, \alpha, \beta \).  We have that \( \refl P = (\refl_{st} B; s, t; \beta, \alpha) \).  But from the assumption, \( S \) is reflectively closed, and therefore \( \refl_{s,t} B \) is some (possibly different) block decomposition in \( S \).  We know from \autoref{lem:reflofsub} that \( \refl P \) defines a subsequence on \( \refl_{s,t} B \).  Therefore \( \refl P \in T \), as required.
		\end{proof}
	\end{Lem}
	
	We know from \autoref{lem:reflinvolution} that \( \refl \) is an involution on \( T \), meaning that \( \refl^2 = \id_T \).  We can consider the group \( G = \Set{\id_T, \refl} \), and its action on the set \( T \) of subsequences.
	
	\begin{Lem}
		The group \( G = \Set{\id_T, \refl} \) acts on \( T \).
		
		\begin{proof}
			This is clear since \( G \) is a group of functions, and the action is function application.  The rule for evaluating \( (f \circ g)(x) \) as \( f(g(x)) \) is one of the condition for a group action.  That the function \( \id_T \) is the identity function on \( T \) is the other condition for a group action.
		\end{proof}
	\end{Lem}
	
	\begin{Lem}\label{lem:subsequenceorbit}
		The set \( T \) breaks up into orbits of size \( \leq 2 \) under the action of \( G = \Set{\id_T, \refl} \).
		
		\begin{proof}
			By the Orbit-Stabilizer theorem, the size of an orbit under this action divides the size of \( G \), which is 2.
		\end{proof}
	\end{Lem}
	
	\begin{Lem}\label{lem:oddsubseq_trivial}
		Let \( \mathcal{O} \) be an orbit of \( T \) under \( G \), which has size 1.  Then the subsequence in \( \mathcal{O} \) is trivial, since it has odd length.  (Recall from \autoref{def:trivialsubsequence}, this means the end points of the subsequence are equal.)
		
		\begin{proof}
			Suppose \( \mathcal{O} = \Set{ P } \).  Then we must have \( \refl P = P \), so the subsequence \( P \) in \( \mathcal{O} \) is a fixed point of \( \refl \).  Now, since \( \mathcal{O} \) has odd length, we know from \autoref{lem:subseqfix_trivial} that it is trivial.
		\end{proof}
	\end{Lem}
	
	\begin{Lem}\label{lem:bothtrivial_or_bothnontrivial}
		Suppose that \( \mathcal{O} \) is an orbit of \( T \) under \( G \), which has size 2.  Then either \( \mathcal{O} \) contains two trivial subsequences, or it contains two non-trivial subsequences.
		
		\begin{proof}
			Suppose that \( \mathcal{O} = \Set{ P_1, P_2 } \), and that \( P_1 = (B; s, t; \alpha, \beta) \) is non-trivial.  We have therefore that the first digit of \( P_1 \), which is \( B_s^\st + \alpha \), and the last digit of \( P_1 \), which is \( B_t^\en - \beta \), are distinct.
			
			Now compute the first and last digit of \( P_2 = (\refl_{st} B; s, t; \beta, \alpha) \).  In the case where \( B_s^\st = B_t^\en \), we get \( (\refl_{st} B)_s^\st = B_t^\en \) and \( (\refl_{st} B)_t^\en = B_t^\st \), so that the first and last digits of \( P_2 \) are \( B_t^\en + \beta = B_t^\en - \beta \pmod*{2} \), and \( B_s^\st - \alpha = B_s^\st + \alpha \pmod*{2} \).  These are the same as those of \( P_1 \), so are still distinct.
			
			In the case where \( B_s^\st \neq B_t^\en \), we find  \( (\refl_{st} B)_s^\st = 1 - B_t^\en \) and \( (\refl_{st} B)_t^\en = 1 - B_t^\st \), so that the first and last digits of \( P_2 \) are \( 1 - B_t^\en + \beta = 1 + B_t^\en - \beta \pmod*{2} \), and \( 1 - B_s^\st - \alpha = 1 + B_s^\st + \alpha \pmod*{2} \).  Since these are the opposite of those of \( P_1 \), they are also distinct.
		\end{proof}
	\end{Lem}
	
	\begin{Lem}\label{lem:subseqnegative}
		Let \( \mathcal{O} \) be an orbit of \( T \) under \( G \), which consists of two non-trivial subsequences.  Then the quotient sequences determined by these subsequences are equal, and the integrals of the subsequences are negatives of each other.
		
		\begin{proof}
			Let the two subsequence be \( P_1 = (B; s, t; \alpha, \beta) \) and \( P_2 = \refl P_1 = (\refl_{st} B; s, t, \beta, \alpha) \).    Say \( B = (\epsilon_1; \; \ell_1, \ldots, \ell_n) \).  Then for \( i < s \) and \( i > t \), the blocks of \( B \) and \( \refl_{s,t} B \) agree, so the quotient sequences agree here.  Since \( P_1 \) is non-trivial, the first and last letters are different.  Suppose \( P \) starts with \( x \), then it ends with \( 1-x \).  Set \( \delta = B_s^\st \).  Then the quotient sequence is
			\[
			W_{\epsilon_1}^{\ell_1} \cdots W_{\epsilon_{s-1}}^{\ell_{s-1}} W_{\delta}^{\alpha + 1} \oplus W_{1-x}^{\beta + 1} W_{\epsilon_{t+1}}^{\ell_{t+1}} \cdots W_{\epsilon_n}^{\ell_n} \, .
			\]
			Since \( W_{\delta}^{\alpha+1} \) ends with \( x \), and \( W_{1-x}^{\beta+1} \) starts with \( 1- x \), we have 
			\[
			W_{\delta}^{\alpha + 1} \oplus W_{1-x}^{\beta + 1} = W_{\delta}^{\alpha+\beta+2} \, .
			\]
			So the blocks in \( B \) are joined by the word \( W_\delta^{\alpha + \beta + 2} \).
			
			But for the same reason, the blocks in \( \refl_{s,t} B \) are joined by \( W_\epsilon^{\beta + \alpha + 2} \), where \( \epsilon = (\refl_{st} B)_s^\st = B_s^\st \), by the definition of \( \refl_{s,t} \).  Therefore the quotient sequences \( Q_1 \) from \( P_1 \) and \( Q_2 \) from \( P_2 \) are both identical.  So we certainly have \( I^\mot(Q_1) = I^\mot(Q_2) \).  \medskip
			
			Using \autoref{lem:revordual}, we know that the subsequences \( P_1 \) and \( P_2 \) are either the reverse, or the dual, of each other.  If \( P_2 \) is the reverse of \( P_1 \), then by the \hyperref[ppty:reversalofpaths]{reversal of paths} property from \autoref{ppty:propertiesofintegrals}, we have \( I^\lmot(P_1) = -I^\lmot(P_2) \) since \( P_1 \) and \( P_2 \) have odd length.  If \( P_2 \) is the dual of \( P_1 \) then, by \hyperref[ppty:duality]{duality}, we also have \( I^\lmot(P_1) = -I^\lmot(P_2) \), since \( P_1 \) and \( P_2 \) have odd length.
		\end{proof}
	\end{Lem}
	
	\begin{Lem}\label{lem:orbitzero}
		Let \( \mathcal{O} \) be an orbit of \( T \) under \( G \).  Then the sum of the terms that \( \mathcal{O} \) gives rise to in \( D_{<N} \) is 0.
		
		\begin{proof}
			If \( \mathcal{O} \) has size 1 then, by \autoref{lem:oddsubseq_trivial}, the subsequence in \( \mathcal{O} \) is trivial, and the orbit \( \mathcal{O} \) contributes \( 0 \) to \( D_{<N} \).\medskip
			
			If \( \mathcal{O} \) has size 2 and the two subsequences it contains are trivial, then the orbit \( \mathcal{O} \) contributes \( 0 \) to \( D_{<N} \).  Otherwise, by \autoref{lem:bothtrivial_or_bothnontrivial}, the two subsequences \( P_1 \) and \( P_2 \) in \( \mathcal{O} \) are non-trivial.  But then by \autoref{lem:subseqnegative} we have \( I^\lmot(P_1) = -I^\lmot(P_2) \), and \( I^\mot(Q_1) = I^\mot(Q_2) \), where \( Q_i \) is the quotient sequence obtained from \( P_i \). Then the orbit \( \mathcal{O} \) contributes
			\begin{align*}
			& I^\lmot(P_1) \otimes I^\mot(Q_1) + I^\lmot(P_2) \otimes I^\mot(Q_2)  \\
			& = I^\lmot(P_1) \otimes I^\mot(Q_1) - I^\lmot(P_1) \otimes I^\mot(Q_1) \\
			& = 0 \qedhere
			\, . \end{align*}
		\end{proof}
	\end{Lem}
	
	We can now present the proof of the following theorem.
	\thmreflectivelyclosed
	
	\begin{proof}
		The goal is to compute \( D_{<N} \) on the sum \( \sum_{s \in S} I_\bl^\mot(s) \).  Since the coefficients of all the integrals in the sum are \( + 1 \), the terms of \( D_{<N} \) arise exactly from the set of all odd subsequences on the block decompositions in \( S \).  Write \( T \) for the set of all odd subsequences on \( S \).
		
		By \autoref{lem:subsequenceorbit} we know that the set \( T \) breaks up into orbits of size \( \leq 2 \) under the action of the group \( \Set{ \id_T, \refl } \) generated by the reflection operator.  From \autoref{lem:orbitzero} we know that each of these orbits contributes \( 0 \) to \( D_{<N} \).
		
		Therefore \( D_{<N} \left( \sum\nolimits_{s \in S} I_\bl^\mot(s) \right) = 0 \), and by Brown's characterisation of \( \ker D_{<N} \), \autoref{thm:kerDN}, we have \( \sum\nolimits_{s \in S} I_\bl^\mot(s) \in \zeta^\mot(N) \Q \), as required.	
	\end{proof}
	
	\begin{Rem}\label{rem:oddweight:trivial}
		In any reflectively closed set of odd weight, the block decomposition \( B \) and the reflection \( \refl_{1,n} B \) give rise to integrals with opposite value, by duality, as in \autoref{rem:dualityr1n}.  Therefore in odd weight, the integrals of a reflectively closed set of block decompositions necessarily sum to 0.  So the above theorem is \emph{trivially} true in odd weight.
	\end{Rem}
	
	From this theorem we naturally obtain a way of \emph{generating} identities on iterated integrals, by forming the reflectively closure of some initial starting set of iterated integrals.
		
		\begin{Cor}\label{cor:generateidentity}
			Let \( \widetilde{S} = \Set{ I^\mot(w_i) } \) be a set of iterated integrals, with corresponding block decompositions \( S = \Set{ B_i } \).  Suppose that \( S \) consists of block decompositions with a fixed (even) weight \( N \), and a fixed number of blocks, but that \( S \) is not necessarily reflectively closed.  Then
			\[
			\sum\nolimits_{s \in \langle S \rangle_{\refl}} I_\bl^\mot(s) \in \zeta^\mot(N) \Q \, .
			\]
		\end{Cor}
		
	And in particular, when this is applied to the a single block decomposition \( B = (\ell_1, \ldots, \ell_n) \) , we obtain the following.  This identity will lend some support to the cyclic insertion conjecture (\autoref{conj:generalcyclicinsertion}), by showing that a sufficiently symmetrised version always holds.
	
	\corsymmetricinsertion
		
	\section{The generalised cyclic insertion conjecture}
	\label{sec:cyclicinsertion}
	
	\subsection{Statement of the conjecture}
	\label{sec:cyclicinsertion:statement}
		
	As shown in \autoref{sec:intro:eg1} and \autoref{sec:intro:eg2}, in the introduction, both Hoffman's identity and the BBBL cyclic insertion conjecture take the form
	\[
		\sum_{\sigma \in C_n} I_\bl(\ell_{\sigma(1)},\ldots,\ell_{\sigma(n)}) = I_\bl(\wt + 2) \, ,
	\]
	for some choice of lengths \( [\ell_i] \).  Numerical experiments with various other choices for \( [\ell_i] \) lead to us the following conjecture.
	
	{
	\printextratrue
	\conjcyclicinsertioneasy
	\printextrafalse
	}

	\begin{Rem}
		The restriction that there is no consecutive pair \( (\ell_i, \ell_{i+1}) = (1,1) \) is necessary, as the following examples show.
		
		For the even weight case \( [\ell_i] = [1,1,2,3,3] \), of weight 8, we find that
		\[
		\left(\frac{\pi^8}{9!}\right)^{-1} \sum_{\sigma \in C_5} I_\bl(\ell_{\sigma(1)}, \ldots, \ell_{\sigma(5)}) = 27.89973142\ldots \overset{?}{\not\in} \Q \, .
		\]
		
		For the odd weight case \( [\ell_i] = [1,1,2,3] \), of weight 7, we find numerically that
		\[
			\sum_{\sigma \in C_4} I_\bl(\ell_{\sigma(1)}, \ldots, \ell_{\sigma(4)}) \eqN 2\zeta(2)\zeta(3) \neq 0 \, .
		\]
		
		Although we cannot say for certain that the result in the weight 8 case is irrational, numerical evaluation is sufficient to prove it is \( \neq 1 \).  Similarly in the odd weight case, the result is proven to be \( \neq 0 \).  We conclude \autoref{conj:generalcyclicinsertion} fails for these choices of \( [\ell_i] \).
	\end{Rem}
	
	With a little more work, these results can eventually be rewritten in the following more suggestive form
	\begin{align*}
	\sum_{C_5} I_\bl(1, 1, 2, 3, 3) &\eqN -2 I_\bl(4) I_\bl(2,3,3) + I_\bl(\wt + 2) \\
	\sum_{C_4} I_\bl(1, 1, 2, 3) &\eqN -2 I_\bl(4) I_\bl(2,3) + I_\bl(\wt + 2) \, .
	\end{align*}
	In both cases, the \( I_\bl(\wt+2) \) term from the basic version of the generalised cyclic insertion conjecture (\autoref{conj:generalcyclicinsertion}) is, in some sense, the leading term in the result.  The `higher order corrections' consist entirely of products terms, moreover they involve subsequences of the original \( [\ell_i] \) obtained by deleting the \( (\ell_1,\ell_2) = (1,1) \) `divergence'.
	
	Further numerical experimentation with other choices for \( [\ell_i] \), including cases where \( [\ell_i] = [1, 1, 1, \ldots] \) contains \( 3 \) consecutive 1's or more, lead us to the following more complete version of the cyclic insertion conjecture.
	
	\begin{Conj}[Generalised cyclic insertion, full version]
		\label{conj:generalcyclicinsertionfull}
		Let \( B = (0; \ell_1,\ldots,\ell_n) \) be a non-trivial block decomposition of weight \( N \).  Let 
		\[
			\mathfrak{L}_k \coloneqq \Set{ [m_{k+1}, \ldots, m_n] | \text{\( [\smash{\underbrace{1, \ldots, 1}_{\text{\( k \) times}}}, m_{k+1}, \ldots, m_n] \) is a cyclic permutation of \( [\ell_1, \ldots, \ell_n] \)}} \, . \mathllap{\phantom{\underbrace{x}_{x}}}
		\]
		The set \( \mathfrak{L}_k \) consists of those cyclic permutations of \( [\ell_1,\ldots,\ell_n] \) which begin with \( k \) consecutive 1's.  We then delete the \( k \) consecutive 1's from the beginning of each element. \medskip
		
		Then the following equality holds
		\[
			\sum_{\sigma \in C_n} I_\bl(\ell_{\sigma(1)}, \ldots, \ell_{\sigma(n)}) \eqN I_\bl(N + 2) - \sum_{k = 1}^{\lfloor n/2 \rfloor} (-1)^{k} \frac{2(2\pi)^{2k}}{(2k+2)!} \sum_{\vec{m} \in \mathfrak{L}_{2k}} I_\bl(\vec{m}) \, .
		\]
	\end{Conj}
	 
	\begin{Rem}
		Curiously the coefficient \( \frac{2(2\pi)^{2k}}{(2k+2)!} \) can also be expressed in the following ways
		\[
		\frac{2(2\pi)^{2k}}{(2k+2)!} = -\frac{\zeta(2k+1)}{(2k+1)(2k+2) \zeta'(-2k)} = \frac{2^{2k+1}}{2k+2} \zeta(\{2\}^{k}) = (-1)^k \frac{2^{2k+1}}{2k+2} I_\bl(2k + 2) \, .
		\]	
		Perhaps these alternative forms suggest some avenues for how one could try to prove this identity exactly.
	\end{Rem}

	\begin{Rem}
		If the list \( [\ell_1, \ldots, \ell_n] \) contains no consecutive pair \( (\ell_i, \ell_{i+1}) = (1,1) \), including \( (\ell_n, \ell_1) \), then \( \mathfrak{L}_{k} = \emptyset \) for all \( k \geq 2 \).  This means all of the product terms in \autoref{conj:generalcyclicinsertionfull} vanish, and \autoref{conj:generalcyclicinsertionfull} does indeed reduce to \autoref{conj:generalcyclicinsertion}.
	\end{Rem}
	
	In \autoref{sec:cyclicinsertion:examples} we will present some examples of both the basic cyclic insertion conjecture (\autoref{conj:generalcyclicinsertion}) and the full cyclic insertion conjecture (\autoref{conj:generalcyclicinsertionfull}).  This will explain how the conjecture works, and the illustrate the range of results it can produce.
	
	In \autoref{sec:123mzv:eg} we will also provide numerous examples for the special subclass of so-called 123-MZV's.  For these 123-MZV's, the conjectural cyclic insertion identities do not need to be shuffle regularised before converting back to MZV's and so retain more structure.
	
	\subsection{Evidence and outlook for the conjecture}
	\label{sec:cyclicinsertion:evidence}
	
		The generalised cyclic insertion conjecture in either incarnation, basic (\autoref{conj:generalcyclicinsertion}) or full (\autoref{conj:generalcyclicinsertionfull}) is unproven.  Even the special cases of the BBBL cyclic insertion conjecture (\autoref{conj:bbblcyclic}) and Hoffman's conjectural identity (\autoref{conj:hoffman}) are not proven.  We can, however, provide various pieces of evidence which lend support to the generalised conjecture.
		
		\subsubsection*{Evidence}
		
		Firstly, a redeeming feature of this conjecture is that it unites two previously separate conjectural identities.  By identifying the salient feature, namely the block decomposition (\autoref{sec:alternatingblock}), which underlies both of these previous conjectures we gain a deeper insight into the problem and open up new potential avenues for a proof.
		
		Secondly, Brown's motivic MZV framework leads to a `symmetric insertion' result in the form of \autoref{cor:symins}.  This result shows that some `sufficiently symmetrised' version of the basic cyclic insertion conjecture always holds, up to a rational.  Moreover, the motivic MZV framework can sometimes be used to prove certain  special cases hold, up to a rational, \emph{without} needing to symmetrise first (for example, \autoref{thm:z2n13312:nosym} below).
		
		Finally, I have numerically checked the full version of the cyclic insertion conjecture to 500 decimal places, for \emph{all} non-trivial block decompositions up to weight 18.  The results were positive.  Moreover, in \autoref{sec:furtherrelations:rank} I present a table containing the rank of these cyclic insertion identities (and some other families of identities) up to weight \( N \leq 13 \) (because of computational constraints).  We see that the rank is always bounded above by the expected rank \( 2^{N-2} - d_N \) of all weight \( N \) MZV relations.
		
		\subsubsection*{Outlook}
		
		The question naturally remains of how we can make further progress towards proving the generalised cyclic insertion conjecture. \medskip
		
		The motivic viewpoint allows us a great insight into the structure of the problem.  However, it is not yet powerful enough to resolve the conjecture entirely by itself.  Suppose that somehow we manage to show that \( D_{<N} \) does vanish on any `basic type' cyclic insertion identity, giving
		\[
			D_{<N} \left( \sum_{\sigma \in C_n} I_\bl(\ell_{\sigma(1)}, \ldots, \ell_{\sigma(n)}) \right) \overset{!}{=} 0 \, .
		\]  Then Brown's characterisation of \( \ker D_{<N} \) in \autoref{thm:kerDN} tells us only that
		\[
			\sum_{\sigma \in C_n} I_\bl(\ell_{\sigma(1)}, \ldots, \ell_{\sigma(n)}) = q \zeta^\mot(\wt) \, ,
		\]
		for some \( q \in \Q \).  
		
		For specific \( [\ell_i] \), we can always numerically evaluate both sides of the resulting identity to determine a highly accurate approximation for \( q \).  But in order to prove the cyclic insertion conjecture exactly, we need to exactly evaluate \( q \) for arbitrary \( [\ell_i] \), and the motivic framework is currently unable to do this.  At some point we would need to rely on an \emph{analytically} proven identity and use this to establish the motivic version.
		
		Indeed, in \cite{glanois2016periods}, Glanois already identifies the same issue.  There Glanois introduces the notion of `families of motivic identities stable under the derivations \( D_{2r+1} \)'.  She uses this to \emph{lift} certain analytically proven families of identities to motivic versions, but she explains the need to already have an analytic version of the identity in order to fix the rational at each step of the motivic lifting.
		
		Nevertheless, this opens up a potential avenue for some kind of `semi-proof' of the generalised cyclic insertion conjecture, up to \( \Q \), \emph{conditional} on an exact version of cyclic insertion at lower weights.  Indeed, we can show that the `basic type' cyclic insertion identities satisfy the following stability under the derivation maps \( D_{2r+1} \).
		
		\begin{Prop}[Stability of cyclic insertion]
			\label{prop:stability}
			The cyclic insertion conjecture is stable under the derivations \( D_{2r+1} \).  More precisely, consider the cyclic insertion identity 
			\[
			Z \coloneqq \sum_{\sigma \in C_n} I_\bl^\mot(\ell_{\sigma(1)}, \ldots, \ell_{\sigma(n)}) \, .
			\]
			The computation of \( D_{2r+1} Z \) results in a sum of terms of the following form.
			\[
			I_\bl^\lmot(a_1,\ldots,a_m) \otimes \sum_{\sigma \in C_k} I_\bl^\mot(b_{\sigma(1)}, \ldots, b_{\sigma(k)}) \, ,
			\]
			where \( m + k = n + 1\), and each \( b_i \) is an \( \ell_i \), or a sum of the form \( \alpha + \beta + 2 \) with \( \alpha \) and \( \beta \) as in \autoref{def:subsequence}.  In particular \( [b_i] \) never has more consecutive 1's than \( [\ell_i] \) in the original identity.
			
			\begin{proof}
				Consider a subsequence marked out on some term of \( Z = \sum_{\sigma \in C_n} I^\mot_\bl(\ell_{\sigma(1)}, \ldots, \ell_{\sigma(n)}) \).  Without loss of generality, we can assume it is marked out on the first term, to obtain the following
				\[
					I^\mot_\bl(\ell_1,\ldots,\underbracket{\overbrace{\ell_i}^{\mathclap{\text{\(\alpha\) from start}}\hspace{0.2cm}},\ldots,\overbrace{\ell_j}^{\mathclap{\hspace{0.2cm}\text{\( \beta \) from end}}}}_{\text{subsequence}},\ldots,\ell_k) \, .
				\]
				Under certain cyclic permutations of the \( [\ell_i] \), the blocks \( \ell_i, \ldots, \ell_j \) are contiguous, and this subsequence has an image there.  These cyclic permutations are the \( n + i-j \) permutations which begin with \( \ell_1, \ell_2, \ldots, \ell_i, \ell_{j+1}, \ldots, \ell_n \).
				
				For each of these permutations, the image subsequence gives rise to the same \( \lmot \)-factor \[ I_\bl^\lmot(\ell_i-\alpha, \ell_{i+1}, \ldots, \ell_{j-1}, \ell_j - \beta) \, . 
				\]
				On the other hand, the resulting quotient sequence is a cyclic permutation of
				\[
					 I_\bl^\mot(\ell_1, \ldots, \ell_{i-1}, \alpha+\beta + 2, \ell_{j+1}, \ldots, \ell_n) \, ,
				\]
				namely the one which starts with the same \( \ell_k \) as the original permutation.  Now group together those terms with this common \( \lmot \)-factor.  This shows that \( D_{2r+1} Z \) consists of a sum of terms of the form 
				\[
					I_\bl^\lmot(a_1,\ldots,a_m) \otimes \sum_{\sigma \in C_k} I_\bl^\mot(b_{\sigma(1)}, \ldots, b_{\sigma(k)}) \, ,
				\]
				as claimed.
				
				Finally, the \( \lmot \)-factor has \( m = j-i+1 \) blocks, whereas the \( \mot \)-factor has \( k = n + i - j \) blocks.  And indeed \( m + k = n + 1 \).
			\end{proof}
		\end{Prop}
		
		When the weight is even, the lower weight cyclic insertion identities which appear in \autoref{prop:stability} are of odd weight.  Therefore we obtain the following corollary, conditional on the `basic version' of the odd weight cyclic insertion conjecture holding exactly.
		
		\begin{Cor}[Even weight cyclic insertion]\label{cor:evenweightconditional}
			Let \( B = (0; \ell_1,\ldots,\ell_n) \) be a non-trivial block decomposition of even weight \( N \), with no pair of consecutive blocks \( (\ell_i, \ell_{i+1}) = (1,1) \).  Suppose that the `basic version' of the odd-weight cyclic insertion conjecture holds for all weights \( <  N - 2 \).  Then
			\[
				\sum_{\sigma \in C_n} I_\bl^\mot(\ell_{\sigma(1)}, \ldots, \ell_{\sigma(n)}) \in I_\bl^\mot(N + 2 ) \Q \, .
			\]
			Therefore the `basic version' even weight cyclic insertion conjecture does indeed hold, up to a rational.
			
			\begin{proof}
				The computation in \autoref{prop:stability} shows that this \( D_{<N} \) is a sum of odd weight cyclic insertion identities of weight \( < N-2 \), in the \( \mot \)-factor.  By assumption each of these vanishes, using \autoref{conj:generalcyclicinsertion}. Therefore
				\[
					D_{<N} \left( \sum_{\sigma \in C_n} I_\bl^\mot(\ell_{\sigma(1)}, \ldots, \ell_{\sigma(n)}) \right) = 0 \, ,
				\]
				and by the characterisation of \( \ker D_{<N} \) in \autoref{thm:kerDN}, the result is a rational multiple of \( \zeta^\mot(N) \).  Or equivalently, it is a rational multiple of \( I^\mot_\bl(N+2) \).
			\end{proof}
		\end{Cor}
		
		\begin{Rem}
			In the case of odd weight \( N \), we do not yet have such an easy corollary.  When computing \( D_{2r+1} \) we obtain a combination of the form \[
				\left( \sum I_\bl^\lmot(a_1,\ldots,a_k)  \right) \otimes I_\bl^\mot(M) \, ,
			\]
			for \( M = N - (2r-1) \) an even integer.  Since \( I^\bl_\mot(M) \neq 0 \),  we have to \emph{explicitly} compute terms modulo products in the \( \lmot \)-factor.
			
			For example the computation of \( D_7 \left( \sum_{C_4} I^\mot_\bl(2, 10, 3, 2) \right) \) produces the following result
			\[
				\big(I_\bl^\lmot(6, 3) + I_\bl^\lmot(3, 3, 2, 1) + I_\bl^\lmot(2, 3, 2, 3) + I_\bl^\mot(1, 2, 2, 4) \big) \otimes I_\bl^\mot(10) \, .
			\]
			To see that \( D_7 \) vanishes, we must recognise that this \( \lmot \)-factor can be rewritten as
			\[
				= -\zeta^\lmot(2)\zeta^\lmot(2,3) - 2\zeta^\lmot(2)\zeta^\lmot(3,2) + 2\zeta^\lmot(3) \zeta^\lmot(2,2) \, .
			\]
			Now \( D_7 \) is seen to vanish because in the \( \lmot \)-factor we work modulo products. \medskip
		\end{Rem}
		
		It is not yet clear to me why the \( \lmot \)-factor \emph{should} always vanish modulo products, when we apply the derivations \( D_{2r+1} \) to a cyclic insertion identity.  However we do have the following observation, about the structure of the \( \lmot \)-factor, thanks to Panzer.
		
		\begin{Obs}[Panzer]\label{obs:dodd:simplification}
			The very special form of \( D_{2r+1} \sum_{\sigma \in C_n} I_\bl^\mot(\ell_{\sigma(1)}, \ldots, \ell_{\sigma(n)}) \) in the odd weight case \( N = 2t + 1 \) does already yield some information about 
			\[
				Z^\mot \coloneqq \sum_{\sigma \in C_n} I_\bl^\mot(\ell_{\sigma(1)}, \ldots, \ell_{\sigma(n)}) \, .
			\]
			Specifically, we have that
			\[
				Z^\mot = \sum_{k = 1}^t \alpha_k \zeta^\mot(2k+1) \zeta^\mot(\{2\}^{t - k})
			\]
			for some rational coefficients \( \alpha_k \in \Q \).  This means that the computation of \( D_{2k+1} Z^\mot \) necessarily simplifies to \( a_k \zeta^\lmot(2k+1) \otimes \zeta^\mot(\{2\})^{t-k} \), and we are only left to determine whether or not \( a_k = 0 \).
			
			\begin{proof}
				The proof of this makes use of the non-canonical isomorphism \( \phi \) of graded Hopf coalgebras, from the motivic MZV's \( \mathcal{H} \) to the \( f \)-alphabet \( \mathcal{U} \coloneqq \Q\langle f_3, f_5, \ldots \rangle \otimes \Q[f_2] \), and how this isomorphism \( \phi \) interacts with the coproduct.  For further details see \cite{brown2012decomposition}.
				
				The isomorphism \( \phi \) can be chosen so that it is `normalised' in depth 1.  This mean \( \phi(\zeta^\mot(N)) = f_N \), for \( N \geq 2 \), where \( f_{2n} = b_n f_2^n \), and \( b_n \) is given by the coefficient \( b_n \) in \( \zeta(2n) = b_n \zeta(2)^n \).  We have that \( \Delta_\dec \circ \phi = \phi \circ \Delta \), where \( \Delta \) is the usual coaction on motivic MZV's and \( \Delta_\dec \) is the deconcatenation coaction on \( \mathcal{U} = \mathcal{U'} \otimes \Q[f_2] \), which lifts \( \Delta_\dec \) on \( \mathcal{U}' = \Q\langle f_3, f_5, \ldots, \rangle \).  Moreover \( \phi \circ D_{2k+1} \) is obtained by taking the \( (2k+1, N - 2k -1) \)-graded piece of \( \Delta_\dec \circ \phi \), and projecting the first component to \( \mathcal{U'} / (\mathcal{U}'_{>0} \cdot \mathcal{U}'_{>0}) \). \medskip
				
				Suppose that \( Y^\mot = \sum \beta_i \zeta^\mot(\vec{s}_i) \) is some given linear combination of MZV's of weight \( N = 2t + 1 \).  Suppose further, that \( D_{2k+1} Y^\mot = \left( \sum_j \gamma_j \zeta^\mot(\vec{s'}_j) \right) \otimes \zeta^\mot(2t - 2k) \), with \( \gamma_i \in \Q \),  for every \( 1 \leq k < t \).  Then Panzer claims that the \( f \)-alphabet decomposition of \( Y^\mot \) must be of the form
				\[
					\phi(Y^\mot) = \sum_{k = 1}^t \beta'_k f_{2k+1} f_2^{t - k} \, ,
				\]
				for some rational numbers \( \beta'_k \in \Q \).
				
				Why is this?  If the \( f \)-alphabet decomposition of \( Y^\mot \) is not of this form, it must contain some term starting with \( \geq 2 \) letters \( f_\text{odd} \), namely \( f_{2k_1+1} f_{2k_2 + 1} \cdots f_{2k_n + 1} f_2^{m} \).  This contributes a term \( f_{2k_1+1} \otimes f_{2k_2+1} \cdots f_{2k_n +1} f_2^m \) to the \( (2k_1+1, N - 2k_1 - 1) \)-graded piece of the deconcatenation coaction.  But no such term exists, since the right hand factor of \( D_{2k_1+1} Y^\mot \) is assumed to only be \( \zeta^\mot(2t - 2k_1) \); this has \( f \)-alphabet decomposition \( b_{t-k_1} f_2^{t-k_1} \), and contains no \( f_\text{odd} \) letters.  So indeed Panzer's claim holds.
				
				Now since \( \phi \) is chosen to be normalised in depth 1, we find
				\[
					Y^\mot = \sum_{k=1}^t \beta'_k \zeta^\mot(2k+1) \zeta^\mot(2)^{t-k} = \sum_{k=1}^t \beta''_k \zeta^\mot(2k+1) \zeta^\mot(\{2\}^{t-k}) \, ,
				\]
				for some \( \beta''_k \in \Q \).  Applying this to 
				\[
					Z^\mot \coloneqq \sum_{\sigma \in C_n} I_\bl^\mot(\ell_{\sigma(1)}, \ldots, \ell_{\sigma(n)})
				\] gives the above result.
			\end{proof}
		\end{Obs}
	
		\subsection{Examples of generalised cyclic insertion}
		\label{sec:cyclicinsertion:examples}
	
		We can easily give examples of cyclic insertion for any block decomposition.  For every block decomposition, we can convert back to MZV's using the regularisation procedure described in \autoref{sec:shufflereg} to get MZV level identities.
		
		\begin{Eg}\label{eg:cyc11224}
			Consider the block lengths \( [\ell_i] = [2, 1, 6, 1, 2] \), which describe the iterated integral \( I(01 \mid 1 \mid 101010 \mid 0 \mid 01) = -\zeta(1, 1, 2, 2, 4) \) of weight 10.  We can apply the basic version of cyclic insertion (\autoref{conj:generalcyclicinsertion}) to the \( [\ell_i] \).  We obtain the following identity on block integrals.
			\begin{align*}
				& I_\bl(2,1,6,1,2) + I_\bl(1,6,1,2,2) + I_\bl(6,1,2,2,1) + {} \\
				& {} + I_\bl(1,2,2,1,6) + I_\bl(2,2,1,6,1) \eqN I_\bl(10)
			\end{align*}
			In terms of iterated integrals, this identity reads
			\begin{align*}
				& I(01 \mid 1 \mid 1010101 \mid 0 \mid 01) + I(0 \mid 010101 \mid 1 \mid 10 \mid 01) + I(010101 \mid 1 \mid 10 \mid 01 \mid 1) + {} \\
				& {} + I(0 \mid 01 \mid 10 \mid 0 \mid 010101) + I(01 \mid 10 \mid 0 \mid 010101 \mid 1) \eqN I(0101010101)
			\end{align*}
			To write this identity in terms of MZV's we need to regularise the last four iterated integrals.  By duality we can combine the second and fifth integral, and regularise to obtain the second line below.  By duality we can also combine the third and fourth integrals, and regularise to obtain the third line below.  This gives the following identity on MZV's.			
			\begin{align*}
			& -\zeta(1,1,2,2,4) + {} \\
			& {} + 2\cdot(-2\zeta(1,4,2,3)-2\zeta(1,4,3,2)  -4\zeta(1,5,2,2)-\zeta(2,4,2,2)) + {} \\
			& {} + 2\cdot(3\zeta(2,2,1,1,4)+\zeta(2,2,1,2,3) + \zeta(2,2,2,1,3) + {} \\
			&	 \quad \quad {} + 2\zeta(2,3,1,1,3)+2\zeta(3,2,1,1,3)) \eqN -\frac{\pi^{10}}{11!} \, .
			\end{align*}
			This identity can be checked using tables of known MZV relations.
		\end{Eg}
	
		Now let us see how the full version of cyclic insertion can be used.
		
		\begin{Eg}
			Consider the block lengths \( [\ell_i] = [1,1,1,2,3] \), of weight 6.  Since \( (\ell_1,\ell_2) = (1,1) \) we must use the full version of the cyclic insertion conjecture to produce an identity.  We first compute the sets \( \mathfrak{L}_2, \mathfrak{L}_4, \ldots \).
			
			To form \( \mathfrak{L}_k \) we select the cyclic permutations of \( [1,1,1,2,3] \) which begin with \( k \) consecutive 1's.  For clarity, we will temporarily write \( \widetilde{\mathfrak{L}_k} \) for the set of these permutations.  We then remove the first \( k \) elements from each permutation in \( \widetilde{\mathfrak{L}_k} \) to obtain \( \mathfrak{L}_k \).  We have
			\begin{align*}
				& \widetilde{\mathfrak{L}_2} = \Set{ [1,1,1,2,3], [1,1,2,3,1] } \text{, and so } \\
				& \mathfrak{L}_2 = \Set{ [1,2,3], [2,3,1] }
			\end{align*}
			Since \( [1,1,1,2,3] \) contains only 3 consecutive 1's, \( \widehat{\mathfrak{L}_4} = \mathfrak{L}_4 = \emptyset \), so we can stop.
			
			Applying \autoref{conj:generalcyclicinsertionfull} give us the following identity on block integrals
			\[
				\sum_{C_5} I_\bl(1,1,1,2,3) \eqN I_\bl(8) +  \frac{2^{3}}{4} I_\bl(4) \left\{ I_\bl(1,2,3) + I_\bl(2,3,1) \right\} \, .
			\]
			Regularising produces the following identity on MZV's.  We have attempted to retain the structure of the block integral identity for comparison purposes.
			\begin{align*}
			& \big( -4 \zeta(1, 5) - 3 \zeta(2, 4) - 2 \zeta(3, 3) - \zeta(4, 2) \big) + {} \\ 
			& {} + \big( 9 \zeta(1, 1, 4) + 6 \zeta(1, 2, 3) + 2 \zeta(1, 3, 2)  + 5 \zeta(2, 1, 3) + 2 \zeta(2, 2, 2) + 2 \zeta(3, 1, 2) \big) + {} \\ 
			& {} + \big( -12 \zeta(1, 5) - 6 \zeta(2, 4) - 4 \zeta(3, 3) - 3 \zeta(4, 2) - \zeta(1, 2, 1, 2) \big) + {} \\
			& {} + ( 10 \zeta(6) ) + ( - \zeta(6) ) \\
			& {} \hspace{2em} \eqN - \zeta(2,2,2) + 2 \zeta(2) \left\{ \big( -2 \zeta(1, 3) - \zeta(2, 2) \big) + ( 3 \zeta(4) ) \right\} \, .
			\end{align*}
			It is possible to simplify this result further, but doing so loses much of the structure of the block integral identity.  The 17 terms on the left hand side combine to 12 terms after simplification.
		\end{Eg}
	
		\begin{Eg}
			Consider the block lengths \( [\ell_i] = [1,1,1,1,2,3] \), of weight 9.  Since \( (\ell_1,\ell_2) = (1,1) \) we must again use the full version of the cyclic insertion conjecture to produce an identity.  We first compute the sets \( \mathfrak{L}_2, \mathfrak{L}_4, \ldots \).
			
			As above, we will write temporarily write \( \widetilde{\mathfrak{L}_k} \) for the set of permutations of \( [\ell_i] \) which begin with \( k \) consecutive 1's.  We then compute \( \mathfrak{L}_k \) from\( \widetilde{\mathfrak{L}_k} \).  We have
			\begin{align*}
				& \widetilde{\mathfrak{L}_2} = \Set{ [1,1,1,1,2,3], [1,1,1,2,3,1], [1,1,2,3,1,1] } \text{, and so } \\
				& \mathfrak{L}_2 = \Set{ [1,1,2,3], [1,2,3,1], [2,3,1,1] } \\[2ex]
				& \widetilde{\mathfrak{L}_4} = \Set{ [1,1,1,1,2,3] } \text{, and so } \\
				& \mathfrak{L}_4 = \Set{ [2,3] }
			\end{align*}
			Since \( [1,1,1,1,2,3] \) contains only 4 consecutive 1's, \( \widetilde{\mathfrak{L}_6} = \mathfrak{L}_6 = \emptyset \), so we can stop.
			
			Applying \autoref{conj:generalcyclicinsertionfull} give us the following identity on block integrals
			\begin{align*}
				\sum_{C_6} I_\bl(1,1,1,1,2,3) \eqN {} & \frac{2^3}{4} I_\bl(4) \left\{ I_\bl(1,1,2,3) + I_\bl(1,2,3,1) + I_\bl(2,3,1,1) \right\} + {} \\
					& \quad + \frac{2^5}{6} I_\bl(6) \left\{ I_\bl(2,3) \right\}
			\end{align*}
			It is possible regularise this identity, and write it in terms of MZV's.  Unfortunately, the regularisation procedure begins to generate a larger and larger number of terms, especially if try to retain the structure of the block integral identity.  Before simplification the regularised left hand side contains 43 terms, even after simplifying it still contains 27 terms.  We therefore leave this identity written in terms of block integrals.
		\end{Eg}
	
		\section{123-MZV's and the generalised cyclic insertion conjecture}
		
		\subsection{Motivation, and properties of 123-MZV's}
		\label{sec:123mzv:cyclicinsertion}
	
		As we see in the examples in \autoref{sec:cyclicinsertion:examples}, it is always possible to write down cyclic insertion identities in terms of MZV's.  However, the regularisation procedure begins to obscure much of the structure of the identity, and begins to generate larger and larger number of terms as the weight increases.    Even restricting to the `basic' version of cycle insertion, the 5 cyclic insertion terms arising from \( [\ell_i] = [2, 1, 6, 1, 2] \) in \autoref{eg:cyc11224} explode to 19 terms, before simplification. 
		
		It would be useful for us to investigate when cyclic insertion identities can be converted directly back to MZV's without regularisation.  These identities will necessarily be shorter, and retain more structure.  This even leads us to a `local' description of cyclic insertion by directly manipulating the arguments of MZV's, without going through the iterated integral representation. \medskip
			
		From \autoref{lem:divergent}, we know that divergent integrals correspond to block decompositions with \( \ell_1 = 1 \), or \( \ell_n = 1 \).  If any \( \ell_i = 1 \), it will eventually be moved into the first position after some cyclic permutation.  So we do not have to regularise if and only if all \( \ell_i > 1 \)
		
		\begin{Def}
			A non-trivial block decomposition \( B = (0; \ell_1,\ldots,\ell_n) \) is called \emph{always convergent} if there all \( \ell_i > 1 \).
		\end{Def}
		
		\begin{Prop}\label{prop:alwaysconvergent:mzv}
			Always convergent block decompositions correspond  directly (i.e. without regularisation) to MZV's \( \zeta(s_1,\ldots,s_k) \) of the following form.
			\begin{enumerate}
				\item Each \( s_i \) is 1, 2, or 3
				\item There is no consecutive pair \( s_i = s_{i+1} = 1 \).
			\end{enumerate}
			
			\begin{proof}
				The iterated integral associated to an always convergent block decomposition \( B \) has no subsequences of the form \( \cdots 0 \mid 0 \mid 0 \cdots \), or \( \cdots 1 \mid 1 \mid 1 \).  The former forbids arguments \( \geq 4 \leftrightarrow 1000 \).  The latter forbids arguments \( 1, 1 \).
				
				Conversely, if all arguments are \( \leq 3 \), and there is no pair \( 1, 1 \), then the subsequences \( 000 \) and \( 111 \) do not appear in the iterated integral.  Therefore all blocks have length \( > 1 \).
			\end{proof}
		\end{Prop}
		
		\begin{Def}[123-MZV]
			An MZV satisfying the conditions in \autoref{prop:alwaysconvergent:mzv} will be called a \emph{123-MZV}.
		\end{Def}
		
		From the block decomposition viewpoint, we can give the following structural description of 123-MZV's.
		
		\begin{Prop}
			\label{prop:123mzvstructure}
			The arguments of a 123-MZV are composed of an arbitrary string formed by concatenating a unique combination of substrings of the following type
			\begin{enumerate}
				\item \( \{2\}^\ell, 3 \), where \( \ell \geq 0 \).  This contributes 1 block with length \( 2 \ell + 3 \).
				\item \( \{2\}^\ell, \underbrace{(1,2), \{2\}^{m_1}, \ldots, (1,2), \{2\}^{m_k}}_{\text{\( k \geq 0 \) repetitions}}, 1, \{2\}^n, 3 \), where \( \ell, m_i, n \geq 0 \).  This contributes \( k + 2 \) blocks with lengths \( 2\ell+2, 2m_1 + 3, \ldots, 2m_k+3, 2n+2 \).
			\end{enumerate}
			Then ending with
			\begin{enumerate}[resume]
				\item \( \{2\}^\ell, \underbrace{(1,2), \{2\}^{m_1}, \ldots, (1,2), \{2\}^{m_k}}_{\text{\( k \geq 0 \) repetitions}} \), where \( \ell, m_i \geq 0 \).  This contributes \( k + 1 \) blocks with lengths \( 2\ell + 2, 2m_1 + 3, \ldots, 2m_k + 3 \).
			\end{enumerate}
			\end{Prop}
			
			Here the notation \( (1,2) \) serves to emphasise that in these MZV's \( (1,2) \) seems to function as one argument.  One should perhaps view it as \( \overline{3} \), the dual of \( 3 \).
			
			\begin{proof}
				In the block decomposition, consider the position of the first block after \( B_1 \) which has \( B_i^\st = 0 \).
				
				Suppose this block occurs at position \( i = 2 \), the first such available position.  Since the first block must end \( 0 \) to make \( B_2^\st = 0 \), the first block must have odd length.  Since we restrict to 123-MZV's, all lengths must be \( > 1 \).  So the word in the integral representation of the MZV must begin
				\[
				0(10)^\ell 10 \mid 0 \cdots = W_0^{2\ell+3} \oplus \cdots \, ,
				\]
				where \( \ell \geq 0 \).  This gives case i).
				
				Now suppose the block occurs at position \( i > 2 \).  Then the first block must end 1, meaning it has even length.  The blocks in position \( j = 2, \ldots, i-2 \) must start 1 and end 1, giving them odd length.  Finally the block in position \( i-1 \) must end 0, giving it even length.  Restricting to 123-MZV's forces all lengths to be \( > 1 \).  So the word in the integral representation of the MZV must begin
				\begin{align*}
				& 0(10)^\ell 1 \mid 10(10)^{m_1}1 \mid \cdots \mid 10(10)^{m_k}1 \mid 10(10)^n \mid 0 \cdots \\ 
				& {} = W_0^{2\ell+2} W_1^{2m_1+3} \cdots W_1^{2m_k+3}  W_1^{2n+2} \oplus \cdots \, ,
				\end{align*}
				where \( \ell, m_1, \ldots, m_k, n \geq 0 \).  This gives case ii).
				
				After dealing with all such blocks, there will be no blocks with \( B_i^\st = 0 \) remaining.  Then the word in the integral representation of the MZV must look as follows
				\begin{align*}
				& 0(10)^\ell 1 \mid 10(10)^{m_1}1 \mid \cdots \mid 10(10)^{m_k}1 \\
				& {} = W_0^{2\ell+2} W_1^{2m_1+3} \cdots W_{2m_k+3} \, ,
				\end{align*}
				where \( \ell, m_i \geq 0 \).
				The first block cannot end with 0, otherwise the second block starts with 0.  Similarly all subsequence blocks must end with 1.  Since they also start with 1, this forces their lengths to be odd.  Finally restricting to 123-MZV's means that all block lengths are \( > 1 \).  This gives case iii) since no further blocks can occur.
			\end{proof}

		With this description of 123-MZV's, it sometimes is convenient to separate out the blocks of \( \{2\}^{\ast} \).  These blocks of 2's are, in some sense, secondary to the structure of the 123-MZV; they can be altered arbitrarily in the ensuing identities, as demonstrated in \autoref{rem:mzv123:blocksof2irrelevant} below.  The main structure of the 123-MZV is in the sequence of arguments \( (1,3) \), \( (1,2) \) and \( 3 \) which appear.
		
		\begin{Not}\label{not:mzv123}
			We shall separate out the blocks of 2's from a 123-MZV, and write
			\[
				\zeta(a_1, \ldots, a_{n-1} \mid b_1, \ldots, b_{n}) \coloneqq \zeta(\{2\}^{b_1}, a_1, \ldots, \{2\}^{b_n}, a_{n-1}, \{2\}^{b_{n}}) \, 
			\]
			where \( a_i \in \Set{1, 3, (1,2)} \).  Notice from \autoref{prop:123mzvstructure}, that \( a_i = (1,2) \), \( a_{i+1} = 3 \) is not a valid combination, and should instead be replaced with \( a_i = 1 \), \( a_{i+1} = 3 \), and a block \( \{2\}^{1} \rightsquigarrow b_{i+1} = 1 \).  Also, \autoref{prop:123mzvstructure} shows that the number \( n = \size{\Set{b_i}} \) counts the number of blocks in the MZV.
		\end{Not}

		We define a `cyclic operator' on 123-MZV's, which acts by manipulating these strings of arguments.  Although the definition seems ad-hoc, it has a more natural interpretation on the corresponding block integrals.
		
		\begin{Def}[Cyclic operator \( \cyc \)]
			Let \( \zeta(s_1, \ldots, s_k) \) be a 123-MZV, as above.  Using the structure of 123-MZV's, we define (via abuse of notation) the operator \( \cyc \) on \( \zeta(s_1,\ldots,s_k) \) as follows.
			\begin{enumerate}
				\item \( \zeta(\{2\}^\ell, 3, \mathrm{rest}) \mapsto -\zeta(\mathrm{rest}, (1,2), \{2\}^l) \).
				\item \( \zeta(\{2\}^\ell, (1,2), \{2\}^{m_1}, \ldots, (1,2), \{2\}^{m_k}, 1, \{2\}^n, 3, \mathrm{rest}) \mapsto \newline \phantom{blah\Big(} (-1)^k \zeta(\mathrm{rest}, 1, \{2\}^\ell, 3, \{2\}^{m_1}, 3, \ldots, \{2\}^{m_k}, 3, \{2\}^n) \)
			\end{enumerate}
			Otherwise, only the final substring appears, and we have
			\begin{enumerate}[resume]
				\item \( \zeta(\{2\}^\ell, (1,2), \{2\}^{m_1}, \ldots, (1,2), \{2\}^{m_k}) \mapsto (-1)^k \zeta(\{2\}^{m_1}, 3, \ldots, \{2\}^{m_k}, 3, \{2\}^\ell) \).
			\end{enumerate}\medskip
			
			In terms of \autoref{not:mzv123}, we have
			\begin{enumerate}
				\item \( \zeta(3, \mathrm{rest} \mid \ell, \mathrm{rest}) \mapsto -\zeta(\mathrm{rest}, (1,2) \mid \mathrm{rest}, \ell) \)
				\item \( \zeta(\{(1,2)\}^k, 1,3, \mathrm{rest} \mid \ell, m_1, \ldots, m_k, n, \mathrm{rest}) \mapsto \newline
				\phantom{blah\Big(} (-1)^k \zeta(\mathrm{rest}, 1, 3, \{3\}^k \mid \mathrm{rest}, \ell, m_1, \ldots, m_k, n) \)
				\item \( \zeta(\{(1,2)\}^k \mid \ell, m_1, \ldots, m_k) \mapsto (-1)^k \zeta(\{3\}^k \mid m_1, \ldots, m_k, \ell) \)
			\end{enumerate}
		\end{Def}
		
		\begin{Rem}[Interpertation of \( \cyc \)]
			From the proof of \autoref{prop:123mzvstructure}, which identifies the structure of 123-MZV's, we obtain the following `interpretation' of \( \cyc \).
			
			Let \( \pm\zeta(s_1,\ldots,s_k) \) be a 123-MZV, whose corresponding iterated integral has block decomposition \( B = (0; \ell_1,\ldots,\ell_n) \).  The sign should be chosen so that this is an equality.  Let \( i > 1 \) be the first position (after \( i = 1 \)), for which \( B_i^\st = 0 \).  Then
			\[
				\cyc \zeta(s_1,\ldots,s_k) = I_\bl(\ell_i, \ell_{i+1}, \ldots, \ell_n, \ell_1, \ldots, \ell_{i-1})
			\]
			is obtained by taking the cyclic permutation of \( [\ell_1, \ldots, \ell_n] \) which starts \( \ell_i \).
			
			Otherwise we have
			\[
				\cyc \zeta(s_1,\ldots,s_k) = I_\bl(\ell_2,\ldots,\ell_n,\ell_1) \, ,
			\]
			if there is no such \( i \).
		\end{Rem}
		
		Suppose we now iterate \( \cyc \) on some 123-MZV \( \pm \zeta(s_1,\ldots,s_k) \) with block decomposition \( B = (\ell_1,\ldots,\ell_n) \).  The sign is again chosen to get equality.  If
		\[
			X = \Set{ i > 1 | B_i^\st = 0 } \text{ , and } Y = \Set{ i > 1 | B_i^\st = 1 } \, ,
		\]
		then we observe the following.  The block corresponding to the \( j \)-th element of \( X \) is moved to position one by applying \( \cyc^j \).  Application of \( \cyc \) changes the starting digit of blocks \( 1, \ldots, i-1 \) when they are moved to the end; the starting digit \( B_1^\st = 0 \) must be flipped to match the last digit \( B_n^\en =  1 \).  Therefore the block corresponding to the \( j \)-th element of \( Y \) is moved to position one by applying \( C^{\size{X} + j} \).
		
		Overall we obtain
		\[
			\sum_{\sigma \in C_n} I_\bl(\ell_{\sigma(1)}, \ldots, \ell_{\sigma(n)}) = \pm \sum_{i=0}^{n-1} \cyc^j \zeta(s_1,\ldots,s_k) \, .
		\]
		
		Therefore, we can recast the cyclic insertion conjecture purely in terms of 123-MZV's and the `local' argument manipulation as follows.
		
		\begin{Conj}[Generalised cyclic insertion for 123-MZV's]\label{conj:cyclic123}
			Let \( \zeta(s_1,\ldots,s_d) \) be a 123-MZV with depth \( d \), and \( n \) blocks (counted using \autoref{prop:123mzvstructure}).  Then
			\[
				\sum_{j=0}^{n-1} \cyc^j \zeta(s_1,\ldots,s_d) \eqN \begin{cases}
					(-1)^{\wt/2 - d}\displaystyle\frac{\pi^\wt}{(\wt+1)!} & \text{if \( \wt \) is even} \\
					0 & \text{if \( \wt \) is odd.} \\
				\end{cases}
			\]
			
			\begin{proof}[Derivation from \autoref{conj:generalcyclicinsertion}]
				The MZV \( (-1)^d \zeta(s_1,\ldots,s_k) \) is equal to some iterated integral \( I(w) \), with coefficient \( +1 \).  This integral has some block decomposition \( B = (0; \ell_1,\ldots,\ell_n) \) with \( n \) blocks.  The previous discussion establishes that
				\[
					(-1)^d \sum_{j=0}^{n-1} \cyc^j \zeta(s_1,\ldots,s_k) = \sum_{\sigma \in C_n} I_\bl(\ell_{\sigma(1)}, \ldots, \ell_{\sigma(n)}) \, .
				\]
				
				We know that 123-MZV's correspond to always convergent block decompositions, so we have all \( \ell_i > 1 \), and we are in the `basic' version of the cyclic insertion conjecture.  The cyclic insertion conjecture (\autoref{conj:generalcyclicinsertion}) evaluates this sum as \( I_\bl(\wt + 2) \).
				
				If the weight is odd, then \( I_\bl(\wt+2) = 0 \) using \autoref{lem:weightandblockstrivial}.  Otherwise,
				\[
					I_\bl(\wt+2) = (-1)^{\wt/2} \zeta(\{2\}^{\wt/2}) = (-1)^{\wt/2} \frac{\pi^\wt}{(\wt+1)!} \, ,
				\]
				since \( \zeta(\{2\}^{\wt/2}) \) has depth \( \wt/2 \).  Rearranging this gives the statement above.
			\end{proof}
		\end{Conj}

		\begin{Rem}\label{rem:mzv123:blocksof2irrelevant}
			If we use \autoref{not:mzv123} to write \( \zeta(s_1,\ldots,s_d) \) as \( \zeta(a_1,\ldots,a_{n-1} \mid b_1,\ldots,b_n) \), then we see the sign in \autoref{conj:cyclic123} depends only on the \( a_i \), and their number \( n - 1 = \size{\Set{a_i}} \).
			
			Indeed, the weight of this MZV is
			\(
				\wt = \sum_{i=1}^{n-1} a_i + 2 \sum_{i=1}^n b_i
			\),	and the depth is \( d = (n - 1) + \sum_{i=1}^n b_i \).  Due to the coefficient 2, the \( b_i \) do not affect the parity of the weight.  If \( \sum_{i=1}^{n-1} a_i \) is odd then so is the weight, and the sign is 0 regardless of the values of the \( b_i \).  Otherwise both \( \sum_{i=1}^{n-1} a_i \) and the weight are even.  Then we find
			\[
				\frac{1}{2}\wt - d = \frac{1}{2} \sum_{i=1}^{n-1} a_i - (n-1) \, ,
			\] which depends only on the \( a_i \) and their number \( n - 1 \), as claimed.
		\end{Rem}
		
		\subsection{Examples of cyclic insertion on 123-MZV's}
		\label{sec:123mzv:eg}
		
		We now show, definitely, how the BBBL cyclic insertion conjecture is a special case of this identity, as is Hoffman's conjectural identity.  This will confirm the signs which appeared in the introduction, in \autoref{sec:intro:eg1} and \autoref{sec:intro:eg2}.
		
		\begin{Eg}[BBBL cyclic insertion]\label{eg:bbblcyclic}
			Consider the MZV 
			\begin{align*}
			z &= \zeta(\{2\}^{b_0}, 1, \{2\}^{b_1}, 3, \ldots, 1, \{2\}^{b_{2n-1}}, 3, \{2\}^{b_{2n}}) \\
			&= \zeta(\{1,3\}^{n} \mid b_0, \ldots, b_{2n}) \, ,
			\end{align*} the first MZV appearing in the BBBL cyclic insertion conjecture.
				
			We compute
			\begin{align*}
				& \cyc \zeta(\{1,3\}^n \mid b_0, \ldots, b_{2n}) = \zeta(\{1,3\}^n \mid b_2, \ldots, b_{2n}, b_0, b_1) \, ,
			\end{align*}
			so that the blocks of 2's are moved cyclically by two steps.  Since there are an odd number of blocks of 2's, stepping around by two steps eventually produces all cyclic permutations of the \( a_i \) when we iterate \( \cyc \).  
			
			We see that \( z \) has depth \( d = 2n + \sum b_i \), and even weight \( \wt = 4n + 2 \sum b_i \).  Therefore the coefficient \( (-1)^{\wt/2 + d} = 1 \).  From \autoref{conj:cyclic123} we obtain
			\[
				\sum_{\sigma \in C_{2n+1}} \zeta(\{1,3\}^{n} \mid b_{\sigma(0)}, \ldots, b_{\sigma(2n)}) \eqN \frac{\pi^\wt}{(\wt+1)!} \, ,
			\]
			which recovers the BBBL cyclic insertion conjecture (\autoref{conj:bbblcyclic}).
			
			Alternatively, we can do this calculation in terms of block decompositions.  We have
			\begin{align*}
				& (-1)^d \zeta(\{1,3\}^n \mid b_0, \ldots, b_{2n}) \\ 
				& \quad = (-1)^d \zeta(\{2\}^{b_{0}}, 1, \{2\}^{b_{1}}, 3, \ldots, 1, \{2\}^{b_{2n-1}}, 3, \{2\}^{b_{2n}}) \\
				& \quad = I_\bl(2b_0+2, 2b_1+2, \ldots, 2b_{2n}+2) \, .
			\end{align*}
			We then sum over cyclic permutations of the blocks, to obtain the same result.
		\end{Eg}
		
		Using \autoref{cor:symins}, we can symmetrise this example to a motivically provable identity.
		
		\begin{Thm}[Symmetrised BBBL, Theorem 3.1 in \cite{charlton2015rational}]
			The following symmetrisation of BBBL cyclic insertion holds.
			\[
				\sum_{\sigma \in S_{2n+1}} \zeta(\{1,3\}^n \mid b_{\sigma(0)}, \ldots, b_{\sigma(2n)}) \eqQone (2n)! \frac{\pi^\wt}{(\wt+1)!} \, .
			\]
		\end{Thm}
			
			Recall from \autoref{sec:notation}, that by \( \eqQone \), we mean the left hand side is proven to be a rational multiple of the right hand side, and this rational is expected to be 1 based on numerical evidence.
			
			\begin{proof}
			According to \autoref{cor:symins}, we obtain a motivically provable identity based on the cyclic insertion expression \( \sum_{\sigma \in C_n} I_\bl(\ell_{\sigma(1)}, \ldots, \ell_{\sigma(n)}) \) by replacing \( \sum_{\sigma \in C_n} \) with \( \sum_{\sigma \in S_n} \).  Doing so with the block decomposition in \autoref{eg:bbblcyclic} shows that the left hand side is a rational multiple \( \pi^\wt \).
			
			The identity is then a sum of \( (2n+1)!/(2n+1) = (2n)! \) cyclic insertion expression, each of which is expected to contribute \( \frac{\pi^\wt}{(\wt+1)!} \) according to \autoref{conj:cyclic123}.  This gives the coefficient above.
			\end{proof}

		\begin{Eg}[Hoffman's identity]
			\label{eg:hoffman33}
			Consider the following MZV
			\begin{align*}
				 z &= \zeta(\{2\}^{b_1}, 3, \{2\}^{b_2}, 3, \{2\}^{b_3}) \\
					&= \zeta(3, 3 \mid b_1, b_2, b_3) \, ,
			\end{align*}
			the `fully generic' version of \( \zeta(3,3,\{2\}^m) \) which appears in Hoffman's conjectural identity.  Iterating \( \cyc \) leads to
			\begin{align*}
			\cyc  z &= -\zeta(3,(1,2) \mid b_2, b_3, b_1) = - \zeta(\{2\}^{b_2}, 3, \{2\}^{b_3}, (1,2), \{2\}^{b_1}) \\
			\cyc^2 z &= \zeta((1,2), (1,2) \mid b_3, b_1, b_2) = \zeta(\{2\}^{b_3}, (1,2), \{2\}^{b_1}, (1,2), \{2\}^{b_2}) \, .
			\end{align*}
			
			The depth of \( z \) is \( d = (b_1 + b_2 + b_3) + 2 \), and the weight \( \wt = 2(b_1 + b_2 + b_3) + 6 \) is even.  This gives the coefficient \( (-1)^{\wt/2 + d} = -1 \).  By \autoref{conj:cyclic123} we obtain
			\begin{align*}
				\zeta(3,3 \mid b_1, b_2, b_3) - \zeta(3,(1,2) \mid b_2,b_3,b_1) + \zeta((1,2),(1,2) \mid b_3,b_1,b_2) \eqN - \frac{\pi^\wt}{(\wt+1)!} \, .
			\end{align*}
			In particular, for \( b_1 = b_2 = 0 \) and \( b_3 = m \), we obtain Hoffman's original conjectural identity (\autoref{conj:hoffman}).
			
			In terms of block decompositions, we have
			\begin{align*}
				& (-1)^d \zeta(3,3 \mid b_1,b_2,b_3) = I_\bl(2b_1+3, 2b_2 + 3, 2b_3+2) \\
				& -(-1)^d \zeta(3,(1,2) \mid b_2,b_3,b_1) = I_\bl(2b_2+3,2b_3+2,2b_1+3) \\
				& (-1)^d \zeta((1,2),(1,2) \mid b_3,b_1,b_2) = I_\bl(2b_3+2, 2a_1+3, 2b_2+3) \, .
			\end{align*}
			By summing these cyclic permutations of the blocks, we obtain the same result.
		\end{Eg}
	
		Symmetrising using \autoref{cor:symins} gives the following.
		
		\begin{Thm}[Symmetrised Hoffman]\label{thm:hoffman33}
			We obtain a motivic proof of Hoffman's identity, up to \( \Q \), by symmetrising and using duality.
			\begin{align*}
				\zeta(3,3 \mid b_1, b_2, b_3) - \zeta(3,(1,2) \mid b_2,b_3,b_1) + \zeta((1,2),(1,2) \mid b_3,b_1,b_2) \eqQone - \frac{\pi^\wt}{(\wt+1)!} \, .
			\end{align*}
			In particular, for \( b_1 = b_2 = 0 \) and \( b_3 = m \) we obtain a motivic proof of Hoffman's conjecture (\autoref{conj:hoffman}), up to a rational.  This confirms Zhao's observation that Hoffman's identity can be motivically proven \cite[p. 188]{zhao2016multiple}.

			\begin{Rem}
				It has since been proven, in Theorem 1 of \cite{hirose2017hoffman}, that the rational in \autoref{thm:hoffman33} above is \( 1 \).  We can therefore replace \( \eqQone \) with \( = \), above.  Moreover, we see the result holds even as an identity on motivic MZV's.
			\end{Rem}
			
			\begin{proof}
				The cyclic version is obtained by applying the cyclic insertion conjecture with the block decomposition \( B = (0; 2b_1+3, 2b_2+3, 2b_3+2) \).  Using \autoref{cor:symins}, we replace \( \sum_{\sigma \in C_3} \) with \( \sum_{\sigma \in S_3} \) to get the following motivically provable identity.
				\begin{align*}
				&\zeta(3, 3 \mid b_1, b_2, b_3)  - \zeta(3, (1,2) \mid b_2, b_3, b_1) + \zeta((1,2),(1,2) \mid b_3, b_1, b_2) + {} \\
				&  {} + \zeta(3, 3, \mid b_2, b_1, b_3) - \zeta(3, (1,2) \mid b_1, b_3, b_2) + \zeta((1,2),(1,2) \mid b_3, b_2, b_1) \\
				 & \quad \eqQone - 2 \frac{\pi^\wt}{(\wt+1)!} \, .
				\end{align*}
				Each row is expected to contribute one lot of \( -\frac{\pi^\wt}{(\wt+1)!} \) by \autoref{conj:cyclic123}, giving the overall coefficient 2.
				
				By duality \( \zeta(3, 3 \mid b_1, b_2, b_3) = \zeta((1,2), (1,2) \mid b_3, b_2, b_1) \), and \( \zeta(3, (1,2) \mid b_2, b_3, b_1) = \zeta(3, (1,2) \mid b_1, b_3, b_2) \), so the first term is equal to the last term, the second term is equal to the fifth, and the third is equal to the fourth.  Combine them, and divide by 2 to obtain the given identity.
			\end{proof}
		\end{Thm}
		
		\begin{Rem}
			Generally, duality shows that any 3 block (even weight) cyclic insertion identity holds up to \( \Q \).  From \autoref{cor:symins}, we know that
			\[
				\sum_{\sigma \in S_3} I_\bl(\ell_{\sigma(1)}, \ell_{\sigma(2)}, \ell_{\sigma(3)}) \in \pi^\wt \Q \, .
			\]
			But by \autoref{rem:dualityr1n}, we know that duality gives \( I_\bl(a,b,c) = I_\bl(c,b,a) \), so that the six terms combine into three pairs.  Therefore we obtain
			\[
				2 \sum_{\sigma \in C_3} I_\bl(\ell_{\sigma(1)}, \ell_{\sigma(2)}, \ell_{\sigma(3)}) \in \pi^\wt \Q \, ,
			\]
			with the same result after dividing by 2.
		\end{Rem}
		
		We now move beyond `known' conjectural identities, and can produce many new conjectural identities, all of which pass extensive numerical testing.  Moreover, we can give a motivically true symmetrisation for each one, using \autoref{cor:symins}.
		
		\begin{Eg}[Generalised Hoffman]\label{eg:generalhoffman}
			The following MZV is a `higher' version of the MZV which generates Hoffman's identity,
			\[
				\zeta(\{3\}^{2n} \mid b_1, \ldots, b_{2n}, c)  = \zeta(\{2\}^{b_1}, 3, \ldots, \{2\}^{b_{2n}}, 3, \{2\}^c) \, .
			\]
			Applying \autoref{conj:cyclic123} to it produces the conjectural result that
			\begin{align*}
				& \sum_{i=1}^{2n+1} (-1)^i \zeta(\{3\}^{2n+1-i}, \{(1,2)\}^{i-1} \mid b_i, \ldots, b_{2n}, c, b_1, \ldots, b_{i-1}) \\
				& \quad \eqN -(-1)^n \frac{\pi^\wt}{(\wt+1)!} \, .
			\end{align*}
		\end{Eg}
		
		Symmetrising produces the following
		
		\begin{Thm}[Symmetrised generalised Hoffman]
			The following identity holds motivically.
			\begin{align*}
			&\sum_{i=1}^{2n+1}  \sum_{\sigma \in S_{2n}}  (-1)^i \zeta(\{3\}^{2n+1-i}, \{(1,2)\}^{i-1} \mid b_{\sigma(i)}, \ldots, b_{\sigma(2n)}, c, b_{\sigma(1)}, \ldots, b_{\sigma(i-1)}) \\
			& \quad \eqQone -(-1)^n (2n)! \frac{\pi^\wt}{(\wt+1)!} \, .
			\end{align*}
			
			\begin{proof}
				The \( i \)-th term in \autoref{eg:generalhoffman} has block decomposition
				\[
					I_\bl(2b_i + 3, \ldots, 2b_{2n}+3, 2c + 2, 2b_{1} + 3, \ldots, 2b_{i-1} + 3)
				\]
				After we sum over all permutation of the block lengths, we can gather those terms where the even length block \( 2c+2 \) is in some fixed position,  \( 2n + 1 - i \) say.  These terms are
				\[
					I_\bl(2b_{\sigma(i)} + 3, \ldots, 2b_{\sigma(2n)}+3, 2c + 2, 2b_{\sigma(1)} + 3, \ldots, 2b_{\sigma(i-1)} + 3) \, .
				\] for any permutation \( \sigma \in S_{2n} \).
				
				For each permutation \( \sigma \in S_{2n} \), the resulting combination is a cyclic insertion expression which contributes \( -(-1)^n \frac{\pi^\wt}{(\wt+1)!} \).  This gives the coefficient \( (2n)! \) in the result.  Alternatively, we count that there are \( (2n+1)! / (2n+1) = (2n)! \) cyclic insertion expressions.
			\end{proof}
		\end{Thm}
		
		We can produce such results for any arbitrary 123-MZV we wish to start with.  If we are willing to regularise, and apply the full version of cyclic insertion (\autoref{conj:generalcyclicinsertionfull}), we can even start with an arbitrary MZV, as in \autoref{eg:cyc11224}.  
		
		\begin{Eg}\label{eg:z1333}
			Consider the MZV
			\[
				\zeta(\{2\}^{b_1}, 1, \{2\}^{b_2}, 3, \{2\}^{b_3}, 3, \{2\}^{b_4}, 3, \{2\}^{b_5}) = \zeta(1, 3, 3, 3, \mid b_1, \ldots, b_5) \, .
			\]
			Iterating \( \cyc \), and using \autoref{conj:cyclic123} leads to the following conjectural identity
			\begin{align*}
			& \zeta(1, 3, 3, 3 \mid b_1, b_2, b_3, b_4, b_5) + \zeta(3, 3, 1, 3 \mid b_3, b_4, b_5, b_1, b_2) + {} \\
			& {} - \zeta(3, 1, 3, (1,2) \mid b_4, b_5, b_1, b_2, b_3) + \zeta(1, 3, (1,2),(1,2) \mid b_5, b_1, b_2, b_3, b_4) + {} \\
			& {} + \zeta((1,2),(1,2), 1, 3 \mid b_2, b_3, b_4, b_5, b_1) \eqN -\frac{\pi^\wt}{(\wt+1)!} \, .
			\end{align*}
		\end{Eg}
		
		\begin{Eg}\label{eg:z13312}
			Consider the MZV
			\[
			\zeta(\{2\}^{b_1}, 1, \{2\}^{b_2}, 3, \{2\}^{b_3}, 3, \{2\}^{b_4}, (1,2) , \{2\}^{b_5}) = \zeta(1 ,3, 3, (1,2) \mid b_1, \ldots, b_5) \, .
			\]
			Iterating \( \cyc \), and using \autoref{conj:cyclic123} leads to the following conjectural identity
			\begin{align*}
				& \zeta(1, 3, 3, (1,2) \mid b_1,b_2,b_3,b_4,b_5) + \zeta(3, (1,2), 1, 3 \mid b_3, b_4, b_5, b_1, b_2) + {} \\
				& - \zeta((1,2), 1, 3, (1,2) \mid b_4, b_5, b_1, b_2, b_3) + \zeta((1,2), 1, 3, 3 \mid b_2, b_3, b_4, b_5, b_1) + {} \\
				& - \zeta(3, 1, 3, 3 \mid b_5, b_1, b_2, b_3, b_4) \eqN \frac{\pi^\wt}{(\wt+1)!}
			\end{align*}
		\end{Eg}
		
		\begin{Not}
			It is convenient to write \( \zeta_{\cyc}(s_1,\ldots,s_k) \) to mean the cyclic insertion expression arising from applying \autoref{conj:cyclic123} to \( \zeta(s_1,\ldots,s_k) \).  Similarly for \( \zeta_{\cyc}(a_1,\ldots,a_n \mid b_1,\ldots,b_{k+1}) \).
		\end{Not}
		
		\begin{Thm}\label{thm:z13312sym}
			The following is a motivically provable symmetrisation of both \autoref{eg:z1333} and \autoref{eg:z13312}.
			\begin{align*}
			& \Sym_{\Set{b_1,b_2, b_5}} \Sym_{\Set{b_3,b_4}} \big( \zeta_{\cyc}(1, 3, 3, 3 \mid b_1, b_2, b_3, b_4, b_5)) \\
			& \hspace{2em} - \zeta_{\cyc}(1, 3, 3, (1,2) \mid b_1, b_2, b_3, b_5, b_4) \big) \eqQone -4! \frac{\pi^\wt}{(\wt+1)!} \, .
			\end{align*}
			
			\begin{proof}
				The MZV \( (-1)^d \zeta(1, 3, 3, 3 \mid b_1, b_2, b_3, b_4, b_5) \) generating the conjectural result in \autoref{eg:z1333} has block decomposition
				\[
					I_\bl(2b_1 + 2, 2b_2 + 2, 2b_3 + 3, 2b_4 + 3, 2b_5 + 2) \, .
				\]
				
				We get the motivically provable symmetrisation by summing over all \( 5! \) permutations of \( b_1, \ldots, b_5 \).  These will be grouped into \( 5! / 5 = 24 \) cyclic insertion identities.  The even lengths can be permuted in \( 3! = 6 \) ways without changing the types of MZV's which can appear.  Similarly, the odd lengths can be permuted in \( 2! = 2 \) ways without changing the types of MZV's which appear.  This reduces the number of distinct permutations to consider to
				\[
					\frac{5!}{5 \cdot 2! \cdot 3!} = 2 \, .
				\]
				These `elementary' permutations are
				\begin{align*}
					& I_\bl(2b_1 + 2, 2b_2 + 2, 2b_3 + 3, 2b_4 + 3, 2b_5 + 2) \text{ and,} \\
					& I_\bl(2b_1 + 2, 2b_2 + 2, 2b_3 + 3, 2b_5 + 2, 2b_4 + 3) \, .
				\end{align*}
				They do not differ by a cyclic shift: in the first the odd length blocks are consecutive, whilst in the second they are separated by 1 block.
				
				When converted back to MZV's, we obtain
				\begin{align*}
					(-1)^d \zeta(1, 3, 3, 3 \mid b_1, b_2, b_3, b_4, b_5) \text{ and, } \\
					-(-1)^d \zeta(1, 3, 3, (1,2) \mid b_1, b_2, b_3, b_5, b_4) \, .
				\end{align*}
				To get all permutations, we then take the cyclic permutations of these `elementary' starting permutations, giving \( \zeta_{\cyc} \)'s, then we symmetrise over \( b_1, b_2, b_5 \), giving \( \Sym_{\Set{b_1,b_2,b_5}} \), and symmetrise over \( b_3, b_4 \), giving \( \Sym_{\Set{b_3,b_4}} \).
				
				From the \autoref{eg:z1333} and \autoref{eg:z13312} we know that each of these \( 24 = 4! \) cyclic insertion expression appearing above is expected to contribute one lot of \( -\frac{\pi^\wt}{(\wt+1)!} \), giving the coefficient above.
			\end{proof}	
		\end{Thm}
	
		\begin{Rem}
			In general when the number of blocks is \emph{composite}, we should take more care when computing the representatives of all permutations of the block lengths modulo cyclic shifts, permutations of even length blocks, and permutations of odd length blocks.  This is because \( S_{\Set{\text{even \( \ell_i \)}}} \times S_{\Set{\text{odd \( \ell_i \)}}} \times C_n \) no longer acts \emph{freely} on \( S_{\Set{\ell_i}} \).
			
			For example, when \( n = 9 \), the following permutation of block lengths
			\[
				(2b_4 + 3, 2b_5 + 3, 2b_6 + 2, 2b_7 + 3, 2b_8 + 3, 2b_9 + 2, 2b_1 + 3, 2b_2 + 3, 2b_3 + 2)
			\]
			is obtained from
			\[
				(2b_1 + 3, 2b_2 + 3, 2b_3 + 2, 2b_4 + 3, 2b_5 + 3, 2b_6 + 2, 2b_7 + 3, 2b_8 + 3, 2b_9 + 2)
			\]
			in two different ways.  We can either cyclically shift left by 3 places, or we can permute the even lengths with \( (b_3, b_6, b_9) \) and the odd lengths with \( (b_1, b_4, b_7)(b_2, b_5, b_6) \).
			
			To work out the number of representatives of \( S_{\Set{\ell_i}} \) modulo \( S_{\Set{\text{even \( \ell_i \)}}} \times S_{\Set{\text{odd \( \ell_i \)}}} \times C_n \) one could use Burnside's counting theorem.  To work out the representatives themselves, one can always start by first quotienting out \( C_n \) and the the larger of \( S_{\Set{\text{even \( \ell_i \)}}} \) or \( S_{\Set{\text{odd \( \ell_i \)}}} \).
		\end{Rem}
		
		Finally, we present an example of the generalised cyclic insertion conjecture in the odd weight case.
		
		\begin{Eg}
			Applying \autoref{conj:cyclic123} to 
			\[
				\zeta(\{2\}^{b_1}, 1, \{2\}^{b_2}, 3, \{2\}^{b_3}, 3, \{2\}^{b_4}) = \zeta(1, 3, 3 \mid b_1, b_2, b_3, b_4)
			\] leads to the following conjectural identity.
			\begin{align*}
			& \zeta(1, 3, 3 \mid b_1, b_2, b_3, b_4) + \zeta(3,1,3 \mid b_3,b_4,b_1,b_2) + {} \\
			& {} - \zeta(1, 3, (1,2) \mid b_4,b_1,b_2,b_3) - \zeta((1,2),1,3 \mid b_2,b_3,b_4,b_1) \eqN 0
			\end{align*}
		\end{Eg}
		
		\begin{Rem}
		With this example, there is no point in producing a symmetrisation.  As explained in \autoref{rem:oddweight:trivial}, all of the terms will cancel pairwise by duality if we try to symmetrise.
		
		Notice, however, that identity does appear in the computation of \( D_{<N} \zeta_{\cyc}(\{1,3\}^2 \mid c_1, \ldots, c_5) \).  This is explained by \autoref{prop:stability}.
		\end{Rem}
		
		In Section 2.6 of \cite{charlton2016identities} I provide a number of further examples of cyclic insertion applied to 123-MZV's, and the of resulting motivically provable symmetrisations.
					
		\subsection{Different symmetrisations of cyclic insertion}
		\label{sec:123mzv:differentsym}
		
		\Autoref{cor:symins} establishes that the following symmetrised version of the cyclic insertion conjecture 
		\[
			\sum_{\sigma \in S_n} I_\bl(\ell_{\sigma(1)}, \ldots, \ell_{\sigma(n)}) \in I_\bl(\wt+2) \Q 
		\] does indeed hold.  However, this type of symmetrisation is not the only one which can be proven motivically.  Sometimes, by good fortune we can motivically prove a cyclic insertion identity on the nose.  Other times, we can obtain different symemtrisations of the corresponding identities.
		
		Currently, I do not have a framework which gives easy proofs of the following results.  They are only established by much more tedious, brute force calculations.  For the full details, I refer to the corresponding results in \cite{charlton2016identities}.
		
		\begin{Thm}[Theorem 2.7.1 in \cite{charlton2016identities}]\label{thm:z2n13312:nosym}
			The following identity, a cyclic insertion identity on the nose, can be motivically proven.  
			\begin{align*}
			& \zeta_{\cyc}(1 ,3, 3, (1,2) \mid m, 0, 0, 0, 0) = \\
			& \quad \zeta(\{2\}^m, 1, 3, 3, (1, 2)) + \zeta(3, (1, 2), 1, \{2\}^m, 3) - \zeta((1,2), 1, \{2\}^m, 3, (1, 2)) + {} \\
			& \quad + \zeta((1, 2), 1, 3, 3, \{2\}^m) - \zeta(3, \{2\}^m, 1, 3, 3) \eqQone \frac{\pi^{2m+10}}{(2m + 11)!}	
			\end{align*}
		\end{Thm}
	
		\begin{Rem}
			As a result of a typo, the statement of Theorem 2.7.1 in \cite{charlton2016identities} refers to the wrong identity.  I mistakenly write \( \zeta_{\cyc}(1,3,3,(1,2) \mid 0, 0, 0, 0, m) \) in the statement of the theorem, even though the proof applies to \( \zeta_{\cyc}(1,3,3,(1,2) \mid m, 0, 0, 0, 0) \).  This is corrected above. 
		\end{Rem}
		
		The identity in \autoref{thm:z2n13312:nosym} would, \emph{a priori}, fall into the symmetrisation given in \autoref{thm:z13312sym} and would involve 6 times as many terms (since some \( b_i = 0 \) are repeated, \( \Sym \) degenerates somewhat).  By good fortune, this is not necessary and the identity already holds motivically.
			
		\begin{proof}
			The proof is straightforward, but tedious: simply compute all of the terms in \( D_{<N} \), and see they cancel pairwise.  The full details of this calculation are presented in Theorem 2.7.1 of \cite{charlton2016identities}.  \medskip
			
			In completing these calculations, we observed some new `features' which might help generalise the symmetrisation framework of \autoref{thm:reflclosedid} and \autoref{cor:symins}.
			
			The first such feature is a type of cancellation by `extended' reflection, as in the following case:
			\begin{align*}
			& I^\mot((01)^{m+1} \mid \overbrace{\underbracket{10 \mid 010} \mid 01}^{\text{reflect}} \mid 101) \text{ cancels with } \\
			& \hspace{8em} \updownarrow \\
			& I^\mot((01)^{m+1} \mid \overbrace{10 \mid \underbracket{010 \mid 01}}^{\text{reflect}} \mid 101) \, ,
			\end{align*}
			when computing \( D_{<N} \) on the term \( \zeta(1,3,3,(1,2) \mid m, 0, 0, 0, 0) = \zeta(\{2\}^m, 1, 3, 3, (1, 2)) \).
			
			The second such feature is a type of cancellation by `splicing', where a substring is cut and pasted to a new location, as in the following case:
			\begin{align*}
				& I^\mot(\underbracket{0\overbrace{10 \mid 0}^{\text{cut}}1} \mid 101 \mid (10)^{m+1} \mid 01) \text{ cancels with } \\
				& \hspace{8em} \updownarrow \\
				& I^\mot(01 \mid 101 \mid (10)^{m+1} \mid \underbracket{0\overbrace{1 \mid 10}^{\mathclap{\text{paste inverse}}}1}) \, ,
			\end{align*}
			when computing \( D_{<N} \) on the terms \( \zeta(3,(1,2),1,3 \mid 0, 0, 0, m, 0) = \zeta(3, (1,2), 1, \{2\}^m, 3) \) and its cyclic shift \( -\zeta((1,2),1,3,(1,2) \mid 0, 0, m, 0, 0) = -\zeta((1,2),1,\{2\}^m,3,(1,2)) \).
		\end{proof}
		
		We can also give a different symmetrisation of the \( \zeta(1,3,3,3 \mid \argdot) \) identity from \autoref{eg:z1333}.  This new symmetrisation is highly reminiscent of the Bowman-Bradley Theorem (\autoref{thm:bowmanbradley}) for \( \zeta(\{1,3\}^n \mid \argdot) \).
		
		\begin{Thm}[Theorem 2.7.13 in \cite{charlton2016identities}]\label{thm:z1333:othersym}
			A different motivically provable symmetrisation for \( \zeta_{\cyc}(1,3,3,3 \mid b_1,\ldots,b_5) \) is given by the following identity.
			\[
			\sum_{\substack{b_1 + \cdots + b_5 = m \\ a_i \geq 0}} \zeta_{\cyc}(1, 3, 3, 3 \mid b_1,b_2,b_3,b_4,b_5) \eqQone -\binom{5+m}{m} \frac{\pi^\wt}{(\wt+1)!}
			\].
			
			\begin{proof}
				This proof is again straightforward, but tedious.  See Theorem 2.7.13 in \cite{charlton2016identities} for the full details.  \medskip
				
				The main step, which forces us to use all compositions \( \sum b_i = m \), rather than just the permutations of some fixed \( b_i \) is the following.
				\begin{align*}
				& I^\mot((01)^{v+1} (10)^{w'} \overbrace{\underbracket{(10)^{w''} (01)^{x+1} \fbox{0}}\fbox{0} (10)^{y+1} (01)}^{\text{reflect}} (01)^{z}) \text{ cancels with } \\
				& \hspace{14em} \updownarrow \\
				& I^\mot((01)^{v+1} (10)^{w'} \overbrace{(10) (01)^{y+1} \fbox{0}\underbracket{\fbox{0} (10)^{x+1} (01)^{w''}}}^{\text{reflect}} (01)^{z}) \, .
				\end{align*}
				This shows that a subsequence in
				\[
				\zeta(1, 3, 3, 3 \mid v, w, x, y, z)
				\] cancels with a subsequence from
				\[
				\zeta(1, 3, 3, 3, \mid v, w', y, x, w''+z-1) \, ,
				\]
				for some \( w' + w'' = w + 1 \).  
				
				These sequences \( (v, w, x, y, z) \) and \( (v, w', y, x, w''+z-1) \) are two different compositions, which cannot be related by a permutation.  However, if we sum over all weak compositions of the parameters \( b_i \), we can guarantee all terms in \( D_{<N} \) cancel.  We find the expected numerical value by counting the number of compositions of \( m \) into 5 non-negative parts; each such composition one cyclic insertion expression to the result.  By a standard result, there are \( \binom{5+m}{m} \) such compositions, each of which is expected to contribute one lot of \( -\frac{\pi^\wt}{(\wt+1)!} \).
			\end{proof}
		\end{Thm}
		
		\begin{Rem}\label{rem:yetmoremotivic}
		This should readily generalise, with the appropriate modifications, to give a motivic proof of
		\begin{align*}
		\sum_{\substack{b_1 + \cdots + b_{2n+1}= m \\ a_i \geq 0}} \zeta_{\cyc}(1, 3, \{3\}^{2n-2} \mid b_1, \ldots, b_{2n+1}) \overset{?}{\eqQone} -\binom{(2n+1)+m}{m} \frac{\pi^\wt}{(\wt+1)!} \, .
		\end{align*}
		
		For the special case \( \zeta_{\cyc}(1,3,3,3 \mid 0, 0, 0, 0, m) \) it \emph{appears} one can give yet more motivically provable symmetrisations.  These symmetrisations all involve a sum of \( \zeta_{\cyc}(1,3,3,3 \mid b_1,\ldots,b_5) \), taken over a very specific list of compositions \( b_1 + \cdots + b_5 = m \).  The simplest example of this is
		\begin{align*}
			& \zeta_{\cyc}(1,3,3,3 \mid 0,0,0,0,m) + \zeta_{\cyc}(1,3,3,3 \mid 0, 0, 0, m, 0) + {} \\
			& \quad {} + \sum_{i=1}^{m-2} \zeta_{\cyc}(1, 3, 3, 3 \mid 0, 0, i, 0, m-i) + {} \label{eqn:further:13332n} \numberthis \\
			& \quad {} + \zeta_{\cyc}(1, 3, 3, 3, \mid 0, 1, 0, m-1, 0) \overset{?}{\eqQone} - (m+1) \frac{\pi^\wt}{(\wt+1)!} \, .
		\end{align*}
		\end{Rem}
	
		\begin{Rem}
			Once the details of a motivic proof of \autoref{eqn:further:13332n} are given, we can try to identify any new cancellation features which arise.  Together with those already identified in the proofs of \autoref{thm:z2n13312:nosym} and \autoref{thm:z1333:othersym}, we should try to use these new cancellation features to further generalise the framework surrounding \autoref{thm:reflclosedid} and \autoref{cor:symins}.  This may allow us to make further progress towards a proof of generalised cyclic insertion conjecture (\autoref{conj:generalcyclicinsertion}) itself.
		\end{Rem}
					
	\section{Other block decomposition relations}
	\label{sec:furtheridentities}
	
	In this final section, we give a number of other numerically verified identities which can be expressed well using the alternating block decomposition.  This is very much in the spirit of the original cyclic insertion conjecture paper \cite{borwein1998combinatorial}.  Indeed we start by recalling another conjectural family of identities put forth that paper \cite{borwein1998combinatorial}. \medskip
	
	In Section 7.2 of \cite{borwein1998combinatorial}, the authors mention that the following identity, distinct from the cyclic insertion conjecture, appears to hold.
	
	\begin{Conj}[Conjecture 2 in \cite{borwein1998combinatorial}]\label{conj:alt13}
		For any non-negative integers \( b_i \geq 0 \), the following identity holds
		\[
			\Alt_{\Set{b_1, b_3, b_5}} \zeta(1, 3, 1, 3 \mid b_1, b_2, b_3, b_4, b_5) \eqN 0 \, .
		\]
	\end{Conj}
	
	The authors of \cite{borwein1998combinatorial} write this conjecture in a different way, which unfortunately obscures the \( \Alt_{\Set{b_i | \text{\( i \) odd}}} \) structure.  They give an expression involving a cyclic shifting of \( \Set{a_1,a_3,a_5} \), which is invariant under \( b_2 \leftrightarrow b_4 \).  Perhaps because of this, they do not give the following generalisation of \autoref{conj:alt13}
	
	\begin{Conj}
		For any non-negative integers \( b_i \geq 0 \), the following identity holds
	\[
		\Alt_{\Set{a_i | \text{\( i \) odd}}} \zeta(\{1,3\}^{2n} \mid b_1, \ldots, b_{2n+1}) \eqN 0 \, .
	\]
	\end{Conj}
	
	Upon writing this conjectural identity in terms of block decompositions, we find that it reads
	\[
		\Alt_{\Set{ b_i | \text{\( i \) odd}}} I_\bl(2b_1 + 2, 2b_2 + 2, \ldots, 2b_{2n+1} + 2) \eqN 0 \, .
	\]
	From here it is natural to consider whether arbitrary block lengths also work.  Indeed they appear to, and so I propose the following conjecture.
	
	\begin{Conj}[Alternate odd positions, `alt-odd']\label{conj:altodd}
		Let \( B = (0; \ell_1,\ldots,\ell_n) \) be any non-trivial block decomposition of even weight \( N \).  Then
		\[
			\Alt_{\Set{\ell_i | \text{\( i \) odd}}} I_\bl^\mot(\ell_1, \ldots, \ell_n) \eqN 0 \, .
		\]
	\end{Conj}

	It appears that there are no restrictions on the block lengths which can appear in \autoref{conj:altodd}.  However, if any \( \ell_{2k+1} = 1 \), we will have to regularise the identity to write it in terms of MZV's.  Instead, if we start with a 123-MZV, where all \( \ell_i > 1 \), we can directly convert the identity back to MZV's and retain the structure of the identity.
	
	\begin{Eg}
		Consider the following 123-MZV which appears in the generalised Hoffman identity (\autoref{eg:generalhoffman}),
		\[ \zeta(\{3\}^{2n-1}, (1,2) \mid b_1, \ldots, b_{2n+1}) \, .
		\]
		It has block decomposition
		\[
			I_\bl(2b_1 + 3, \ldots, 2b_{2n-1} + 3, 2b_{2n} + 2, 2b_{2n+1}+3) \, .
		\]
		consisting of \( 2n + 1 \) blocks.  Since the even length block \( 2 b_{2n} + 2 \) appears in the even index position \( 2n \), it stays fixed when alternating the odd position blocks.  Therefore the MZV type does not change when we apply \autoref{conj:altodd}.  We obtain the conjectural identity
		\[
			\Alt_{\Set{b_i | \text{\( i \) odd}}} \zeta(\{3\}^{2n-1}, (1,2) \mid b_1, \ldots, b_{2n+1}) \eqN 0 \, .
		\]
		
		An analogous result holds for any \( \zeta(\{3\}^{2n-k}, \{(1,2)\}^{k} \mid \argdot) \), as long as the number \( k \) of \( (1,2) \)'s is odd.  This is because in the block decomposition, the even length block \( 2b_{2n - k + 1} + 2 \) occurs in an even index position.  When the odd position blocks are alternated, the even length block remains fixed and the MZV type does not change.  For \( k \) odd we obtain
		\[
			\Alt_{\Set{b_i | \text{\( i \) odd}}} \zeta(\{3\}^{2n-k}, \{(1,2)\}^k \mid b_1, \ldots, b_{2n+1}) \eqN 0 \, .
		\]
	\end{Eg}

	\begin{Rem}
		Indeed, with \( [b_i] = [1, 2, 3, 4, 5, 6, 7, 8, 9] \), the identity
		\begin{align*}
			\Alt_{\Set{b_i | \text{\( i \) odd}}} \zeta(\{3\}^5, \{(1,2)\}^3 \mid b_1, \ldots, b_9) \eqN 0 \, ,
		\end{align*}
		holds numerically to \( > 1000 \) decimal places.  For reference, the individual MVZ's in the identity are each on the order of \( 10^{-126} \), so they conspire to cancel the remaining \( >800 \) decimal places to 0.
	\end{Rem}
		
	\begin{Eg}
		We can also apply the `alt-odd' conjecture to \( \zeta(\{3\}^{2n-k}, \{(1,2)\}^{k} \mid \argdot) \), for even \( k \).  Since the even length block \( 2a_{2n - k + 1} + 2 \) now appears in an odd index position, the resulting identity contains more than one type of MZV.  For example, we obtain
		\begin{align*}
			& \Alt_{\Set{b_1,b_3}} \Big( \zeta(3, 3, 3, 3 \mid b_1, b_2, b_3, b_4, b_5) + \zeta(3, 3, (1,2),(1,2) \mid b_3, b_2, b_5, b_4, b_1) + {} \\
			& \hspace{4em} + \zeta((1,2),(1,2),(1,2),(1,2) \mid b_5, b_2, b_1, b_4, b_3) \Big) \eqN 0 \, ,
		\end{align*}
		by applying this conjecture to
		\[
			\zeta(3,3,3,3 \mid b_1, b_2, b_3, b_4, b_5) = I_\bl(2b_1 + 3, 2b_2 + 3, 2b_3 + 3, 2b_4 + 3, 2b_5 + 2) \, .
		\]
	\end{Eg}

	\subsubsection*{Odd weight `alt-odd'} As yet, there does not appear to be an `easy' odd weight analogue of \autoref{conj:altodd}.  Certainly a number of identities can be found in the odd weight case.  For example, the following identities appear to hold (subject to some restrictions on the allowable block lengths \( \ell_i \), which need to be fully investigated).
	\begin{align}
		& \Alt_{\Set{\ell_1,\ell_3}} \Alt_{\Set{\ell_2,\ell_4}} I_\bl(\ell_1,\ell_2,\ell_3,\ell_4) \eqN 0  \label{eqn:altodd:4block} \\
		& \Alt_{\Set{\ell_1,\ell_4,\ell_6}} \Alt_{\Set{\ell_2,\ell_3,\ell_5}} I_\bl(\ell_1,\ell_2,\ell_3,\ell_4,\ell_5,\ell_6) \eqN 0 \label{eqn:altodd:6block}
	\end{align}
	
	\begin{Rem}
		The 4-block identity, \autoref{eqn:altodd:4block}, tends to fail when consecutive \( (\ell_i, \ell_{i+1}) = (1,1) \) in some term.  The 6-block identity, \autoref{eqn:altodd:6block} appears to be more robust and so far holds for any choice of block lengths.  Unfortunately it is not clear how to generalise these two identities to 8-blocks and beyond.
	\end{Rem}
	
	An alternative candidate for an odd weight analogue of \autoref{conj:altodd} arises when trying to investigate a motivic proof of \autoref{conj:altodd}.  Computation of \( D_{<N} \) leads us to the following identities (subject to some restrictions on the allowable block lengths, which need to be fully investigated).
	\begin{align*}
		& \Alt_{\Set{\ell_1,\ell_3}} (I_\bl(x-\ell_4,\ell_1,\ell_4,\ell_3) + I_\bl(\ell_1,x-\ell_4,\ell_4,\ell_3)) + \\
		& {} - \Alt_{\Set{\ell_1,\ell_3}} (I_\bl(\ell_1,\ell_2,x-\ell_2,\ell_3) + I_\bl(\ell_1,\ell_4,\ell_3,x-\ell_2)) \eqN 0 \\[2ex]
		& \Alt_{\Set{\ell_1,\ell_3,\ell_5}} (I_\bl(x-\ell_4-\ell_6,\ell_1,\ell_4,\ell_3,\ell_6,\ell_5) + I_\bl(\ell_1,x-\ell_4-\ell_6,\ell_4,\ell_3,\ell_6,\ell_5)) + \\
		& - {} \Alt_{\Set{\ell_1,\ell_3,\ell_5}} (I_\bl(\ell_1,\ell_2,x-\ell_2-\ell_6,\ell_3,\ell_6,\ell_5) + I_\bl(\ell_1,\ell_2,\ell_3,x-\ell_2-\ell_6,\ell_6,\ell_5)) + \\
		& + {} \Alt_{\Set{\ell_1,\ell_3,\ell_5}} (I_\bl(\ell_1,\ell_2,\ell_3,\ell_4,x-\ell_2-\ell_4,\ell_5) + I_\bl(\ell_1,\ell_2,\ell_3,\ell_4,\ell_5,x-\ell_2-\ell_4)) \eqN 0 \, .
	\end{align*}
	
	\begin{Rem}
		If we require that the blocks lengths in \emph{every} term are all are \( > 1 \) (both the \( \ell_i > 1 \) and the \( x - \text{sum of \( \ell_i \)'s} > 1 \)), then the identities appear to always hold.  This is a much to strong restriction, as the identities hold in many other cases too.  Unfortunately, the additional parameter \( x \) makes it more difficult to analyse the identities, and determine the true restrictions.
		
		Fortunately, this family does appear to generalise, and to include the previous as a consequence.
	\end{Rem}
	
	The structure of this generalisation needs to be explained carefully.
	
	\subsubsection*{Potential odd weight `alt-odd' candidate}
	
		Let \( \Set{\ell_1,\ldots,\ell_{2n}} \) be some given set of block lengths.  Consider the odd index positions \( \mathcal{O} = \Set{\ell_1,\ell_3,\ldots} \) and the even index positions \( \mathcal{E} = \Set{\ell_2,\ell_4,\ldots} \).  Let \( \mathcal{E}_i \coloneqq \Set{\ell_2, \ell_4,\ldots,\widehat{\ell_{2i}}, \ldots, \ell_{2n}} \) be obtained by dropping the \( i \)-th even position.  Let \( x \) be an additional parameter, with \( x + \sum_{\ell_{2j+1} \in \mathcal{O}} \ell_{2j+1} \) odd.  Finally assume that \( x - \sum_{\ell_{2j} \in \mathcal{E}_i} \ell_{2j} > 0 \) for every \( i \). \medskip
	
		The \( i \)-th row \( \mathcal{R}_i \) of the `alt-odd' identity is obtained as follows.
		\begin{enumerate}[label=Step \arabic*),leftmargin=2cm]
			\item Interleave \( \mathcal{O} \) and \( \mathcal{E}_i \), starting with \( \mathcal{O} \) to obtain the argument string \( \mathcal{A}_i \).
			\item Insert \( x - \sum_{\ell_{2j} \in \mathcal{E}_i} \ell_{2j} \) to the left of \( \ell_{2i-1} \) in \( \mathcal{A}_i \), giving argument string \( \mathcal{B}_i \).
			\item Insert \( x - \sum_{\ell_{2j} \in \mathcal{E}_i} \ell_{2j} \) to the right of \( \ell_{2i-1} \) in \( \mathcal{A}_i \), giving argument string \( \mathcal{C}_i \).
			\item The \( i \)-th row is \( \mathcal{R}_i \coloneqq \Alt_{\mathcal{O}} \left( I_\bl(\mathcal{B}_i) + I_\bl(\mathcal{C}_i) \right) \).
		\end{enumerate}
	
		The odd weight `alt-odd' candidate identity is as follows.
	
		\begin{Conj}[Odd weight `alt-odd' candidate]\label{conj:altodd:oddwtcandidate}
		With the notation as above, the following identity holds (subject to some restrictions on the \( \ell_i \) and \( x \)).
		\[
			\sum_{i=1}^{2n} (-1)^i \mathcal{R}_i \eqN 0 \, .
		\]
		The block lengths in each term of this identity sum to \( x + \sum_{\ell_{2j+1} \in \mathcal{O}} \ell_{2j+1} \), giving us an odd weight combination.
		\end{Conj}
	
	\subsection{Rank of block decomposition identities}
	\label{sec:furtherrelations:rank}
	
	In this final section, we considering what fraction of all MZV relations we obtain from the above block decomposition identities.  We use the full version of the cyclic insertion conjecture \autoref{conj:generalcyclicinsertionfull}, the even weight `alt-odd' conjecture \autoref{conj:altodd}, the odd weight `alt-odd' candidate \autoref{conj:altodd:oddwtcandidate} and the duality relations (as expressed using block decompositions) \autoref{rem:dualityr1n}.
	
	We attempted to prune any duplicate relations cyclic insertion relations, and any duplicate (or trivial) even weight `alt-odd' relations the before computing the rank.  Due to the additional parameter \( x \), it is less obvious how to identify duplicate relations for odd weight `alt-odd' identity, so we did not attempt this.  This explains why the `alt-odd' identity has a much higher initial number of relations at odd weight, in the table below.
	
	The following steps were used to produce the table below.
	\begin{enumerate}[leftmargin=2cm]
		\item[Step 0)] Fix some weight \( N \).  If \( N \) is even, determine all compositions \( \mathcal{X}_N \) of \( N \) into an odd number of parts.  If \( N \) is odd, determine all compositions \( \mathcal{X}_N \) of \( N \) into an even number of parts.
		\item[Step 1)] Select representatives of \( \mathcal{X}_N \) modulo cyclic shifts, apply the full cyclic insertion conjecture (\autoref{conj:generalcyclicinsertionfull}) to each composition, and regularise.  Use the stuffle product of MZV's to resolve the product \( \zeta(n) \ast \zeta(s_1, \ldots, s_k) \) as a linear combination of MZV's.  Compute the rank.
		\item[Step 2a)] For even weight: select representatives of \( \mathcal{X}_N \) modulo alternating the odd positions.  Remove any compositions where the odd positions contain a duplicate value.  Apply the even weight `alt-odd' conjecture (\autoref{conj:altodd}), and regularise.  Compute the rank.
		\item[Step 2b)] For odd weight: compute all compositions \( \mathcal{X}_{\leq N} \) of \( 1, \ldots, N \) into an even number of parts.  Apply the odd weight `alt-odd' candidate (\autoref{conj:altodd:oddwtcandidate}), choosing the additional parameter \( x \) to make the total weight \( N \).  Regularise.  Compute the rank.
		\item[Step 3)] Apply duality to all MZV's of weight \( N \).  Compute the rank.
		\item[Step 4)] Finally compute the overall rank of these 3 families of relations.
	\end{enumerate}
	
	We obtain the following table.
	
	\begin{center}
	\begin{tabular}{c|cc|cc|cc|c|c}
		\multirow{2}{*}{weight} & \multicolumn{6}{c|}{initial number \,\,\, and \,\,\, rank} & \multirow{2}{1.5cm}{overall rank} & \multirow{2}{2.5cm}{expected rank \( 2^{\wt-2} - d_{\wt} \)} \\
		 & \multicolumn{2}{c|}{cyclic} & \multicolumn{2}{c|}{`alt odd'} & \multicolumn{2}{c|}{duality} & & \\ \hline
		4 & 5 & 3 & 1 & 1 & 2 & 1 & 3 & 3 \\
		5 & 7 & 5 & 12 & 4 & 8 & 4 & 5 & 6 \\
		6 & 15 & 13 & 5 & 5 & 12 & 6 & 13 & 14 \\
		7 & 25 & 24 & 72 & 17 & 32 & 16 & 26 & 29 \\
		8 & 51 & 50 & 16 & 16 & 56 & 28 &  56 & 60 \\
		9 & 87 & 86 & 266 & 51 & 128 & 64 & 114 & 123 \\
		10 & 171 & 170 & 41 & 41 & 240 & 120 & 232 & 249 \\
		11 & 306 & 304 & 886 & 142 & 512 & 256  & 465 & 503 \\
		12 & 585 & 584 & 93 & 93 & 992 & 496 & 943 & 1012 \\
		13 & 1081 & 1076 & 2687 & 369 & 2048 & 1024 & \( \ast \) & 2032
	\end{tabular}
	\end{center}

	\begin{Rem}
		Computation of the overall rank at weight 13 failed; the program ran out of memory.  We need to obtain a better understanding of the block length restrictions in the odd weight `alt-odd' candidate, and we need to determine how to identify redundant identities \emph{before} computing the rank.  Once we accomplish this, we will attempt to rerun this calculation and extend it to higher weight.
	\end{Rem}

	Further work to identify more general (alternating) block decomposition identities might help close the relatively small gap between the actual number of relations, and the expected number of relations.  If this gap can be reduced to 0, then we can present yet another potentially complete family of relations on MZV's complementing the many families which are already known.

	{
	\emergencystretch5em
	\printbibliography
	}
	
\end{document}